\documentclass[12pt,letterpaper]{article}

\usepackage{amssymb,latexsym,amsmath,mathrsfs,csquotes}
\usepackage{enumerate}
\usepackage[T1]{fontenc}
\usepackage[dvips,pdftex]{graphicx}
\usepackage{color}
\usepackage{verbatim}
\usepackage[colorlinks=true, urlcolor=blue, linkcolor=blue, citecolor=blue]{hyperref}
\usepackage[left=2.5cm, top=2.5cm, bottom=2.5cm, right=2.5cm]{geometry}

\newtheorem{theorem}{Theorem}[section]
\newtheorem{lemma}[theorem]{Lemma}
\newtheorem{definition}[theorem]{Definition}
\newtheorem{remark}[theorem]{Remark}
\newtheorem{example}[theorem]{Example}

\newtheorem{proposition}[theorem]{Proposition}
\numberwithin{equation}{section}

\def\eqdist{\;{\stackrel{d}{=}}\;}

\definecolor{Red}{rgb}{1,0,0}

\definecolor{Blue}{rgb}{0,0,1}

\title{A Central Limit Theorem for the Optimal Alignments Score in Multiple Random Words}

\author{Ruoting Gong\thanks{Department of Applied Mathematics, Illinois Institute of Technology,
Chicago, Illinois 60616. {\tt Email:\,\url{rgong2@iit.edu}}.} \and Christian Houdr\'{e}\thanks{School of Mathematics, Georgia Institute of Technology, Atlanta, Georgia 30332-0160. {\tt Email:\,\url{houdre@math.gatech.edu}}. Research supported in part by the grant \# 246283 from the Simons Foundation.} \and \"{U}mit I\c{s}lak\thanks{School of Mathematics, Georgia Institute of Technology, Atlanta, Georgia 30332-0160. {\tt Email:\,\url{uislak6@math.gatech.edu}}.}}

\date{March 14, 2016}

\begin{document}

\maketitle

\begin{abstract}
Let $\mathbf{X}^{(1)}_{n},\ldots,\mathbf{X}^{(m)}_{n}$, where $\mathbf{X}^{(i)}_{n}=(X^{(i)}_{1},\ldots,X^{(i)}_{n})$, $i=1,\ldots,m$, be $m$ independent sequences of independent and identically distributed random variables taking their values in a finite alphabet $\mathcal{A}$. Let the score function $S$, defined on $\mathcal{A}^{m}$, be non-negative, bounded, permutation-invariant, and satisfy a bounded differences condition. Under a variance lower-bound assumption, a central limit theorem is proved for the optimal alignments score of the $m$ random words.
\end{abstract}

\vspace{0.3 cm}


\vspace{0.3 cm}


\section{Introduction}

Let $\mathbf{X}^{(j)}_{n}=(X^{(j)}_{1},\ldots,X^{(j)}_{n})$, $j=1,\ldots,m$, be $m$ independent sequences of  i.i.d. random variables drawn from a finite alphabet $\mathcal{A}\subset\mathbb{R}$. Our interest is in studying the similarity/dissimilarity of the random words $(\mathbf{X}^{(i)}_{n})_{i=1}^{m}$, which is central to many areas of applications including computational molecular biology and computational linguistics (cf.~\cite{Capocelli:1990}, \cite{ChristianiniHahn:2007}, \cite{DurbinEddyKroghMitchison:1998}, \cite{LinOch:2004}, \cite{Melamed:1999}, \cite{Pevzner:2000}, \cite{RobinRodopheSchbath:2005}, \cite{SankoffKruskal:1999}, \cite{Waterman:1995}, \cite{YangLi:2003}). For this purpose, let $S:\mathcal{A}^{m}\rightarrow\mathbb{R}^{+}$ be a non-trivial permutation-invariant \emph{score function}, such that
\begin{align}\label{eq:BdScore}
s^{*}:=\sup_{(x_{1},\ldots,x_{m})\in\mathcal{A}^{m}}S(x_{1},\cdots,x_{m})<\infty.
\end{align}
Throughout, the score function $S$ is further assumed to satisfy a bounded differences condition, i.e., let there exist a universal constant $D>0$ such that, for any two elements $(x_{1},\ldots,x_{m})$ and $(y_{1},\ldots,y_{m})\in\mathcal{A}^{m}$ which differ in at most one coordinate,
\begin{align}\label{eq:BddCondScore}
\left|S(x_{1},\cdots,x_{m})-S(y_{1},\cdots,y_{m})\right|\leq D.
\end{align}

Next, define an \emph{alignment} to be an $m$ vectors $(\boldsymbol{\pi}^{(1)},\ldots,\boldsymbol{\pi}^{(m)})$, where each of its coordinates has the same length $k$, for some positive integer $k\le n$, with $\boldsymbol{\pi}^{(j)}=(\pi_{1}^{(j)},\ldots,\pi_{k}^{(j)})$ such that $1\leq\pi^{(j)}_{1}<\pi^{(j)}_{2}<\cdots<\pi^{(j)}_{k}\leq n$ for each $j=1,\ldots,m$. Given a permutation-invariant score function $S:\mathcal{A}^{m}\rightarrow\mathbb{R}^{+}$ satisfying \eqref{eq:BdScore} and \eqref{eq:BddCondScore}, the score of the alignment $(\boldsymbol{\pi}^{(1)},\ldots,\boldsymbol{\pi}^{(m)})$, which aligns $\mathbf{X}_n^{(1)},\ldots,\mathbf{X}_n^{(m)}$, is then defined as
\begin{align*}
U_{(\boldsymbol{\pi}^{(1)},\ldots,\boldsymbol{\pi}^{(m)})}\left(\mathbf{X}_{n}^{(1)};\,\cdots;\mathbf{X}_{n}^{(m)}\right):=\sum_{i=1}^{k}S\left(X^{(1)}_{\pi^{(1)}_{i}},\cdots,X^{(m)}_{\pi^{(m)}_{i}}\right).
\end{align*}
In turn, the optimal score of $\mathbf{X}_n^{(1)},\ldots,\mathbf{X}_n^{(m)}$ is then defined as
\begin{align}\label{eq:OptScore}
L_{n}=L_{n}\left(\mathbf{X}_{n}^{(1)};\,\cdots;\mathbf{X}_{n}^{(m)}\right):=\max_{(\boldsymbol{\pi}^{(1)},\ldots,\boldsymbol{\pi}^{(m)})}\,U_{(\boldsymbol{\pi}^{(1)},\ldots,\boldsymbol{\pi}^{(m)})}\left(\mathbf{X}_{n}^{(1)};\,\cdots;\mathbf{X}_{n}^{(m)}\right),
\end{align}
where the maximum is taken over all the possible alignments. By \eqref{eq:BddCondScore}, if $(\mathbf{x}_{n}^{(i)})_{i=1}^{m}$ and $(\mathbf{y}_{n}^{(i)})_{i=1}^{m}$ are two $m$-tuples of vectors of length $n$ which differ in at most one coordinate, with $\mathbf{x}_{n}^{(i)},\,\mathbf{y}_{n}^{(i)}\in\mathcal{A}^{n}$ for each $i=1,\ldots,m$, then
\begin{align}\label{eq:bddconditionK}
\left|L_{n}\left(\mathbf{x}_{n}^{(1)};\,\cdots;\mathbf{x}_{n}^{(m)}\right)-L_{n}\left(\mathbf{y}_{n}^{(1)};\,\cdots;\mathbf{y}_{n}^{(m)}\right)\right|\leq D.
\end{align}
A particular important example satisfying \eqref{eq:bddconditionK}, with $D=1$, as well as all the other conditions occurs by choosing $S(x_{1},\cdots,x_{m})={\bf 1}_{\{x_{1}=\cdots=x_{m}\}}$, in which case the optimal score is the length of longest common subsequences (LCSs) of the $m$ sequences.
\begin{remark}\label{rem:GapMultiSeqs}
In general, one can also incorporate gap penalties (of arbitrary sign) in our framework by enlarging the alphabet to $\mathcal{A}^{*}:=\mathcal{A}\cup\{\sqcup\}$, where $\sqcup$ is the symbol for a gap, and extending the score function $S$ to a bounded function $S^{*}:(\mathcal{A}^{*})^{m}\rightarrow\mathbb{R}$, with $S^{*}\big|_{\mathcal{A}^{m}}=S$ non-negative, while preserving the permutation-invariant property as well as \eqref{eq:BddCondScore}. The analysis of this extended framework will be similar to ours, except that one needs to modify the definition of alignments when $m\geq 3$. Indeed, when $m=2$, an unaligned letter will be automatically aligned with a gap, so it is sufficient to specify the pairs of letters in an alignment. This is not the case when $m\geq 3$, since a letter can be aligned with different numbers of letters and gaps, which may result in different scores.
\end{remark}

Returning to the probabilistic setting, by a classical super-additivity argument, one easily obtains that
\begin{align}\label{eq:LimitMeanLn}
\lim_{n\rightarrow\infty}\frac{\mathbb{E}(L_{n})}{n}=\gamma^{*}:=\sup_{n\geq 1}\frac{\mathbb{E}(L_{n})}{n},
\end{align}
where $\gamma^{*}\in[0,s^{*}]$ depends on the size of the alphabet, the distribution of $X_{1}^{(1)}$, as well as the score function $S$. We claim that $\gamma^{*}>0$. Indeed, since $S$ is non-trivial, there exists at least one $(x_{1},\ldots,x_{m})\in\mathcal{A}^{m}$ such that
$S(x_{1},\cdots,x_{m})>0$, and therefore
\begin{align*}
\lim_{n\rightarrow\infty}\frac{\mathbb{E}(L_{n})}{n}&\geq\lim_{n\rightarrow\infty}\frac{S(x_{1},\cdots,x_{m})}{n}\sum_{i=1}^{n}\mathbb{P}\left(X_{i}^{(1)}=x_{1},\ldots,X_{i}^{(m)}=x_{m}\right)\\
&=S(x_{1},\cdots,x_{m})\mathbb{P}\left(X_{1}^{(1)}=x_{1},\ldots,X_{1}^{(m)}=x_{m}\right)>0.
\end{align*}

It is the purpose of the present paper to go beyond convergence of the first moment by proving a central limit theorem for $L_{n}$. The methodology of proof will espouse a series of three key steps mirrored in the organization of the paper. To start with, the next section establishes a convergence rate for the first moment of $L_{n}$. In Section \ref{sec:closenesstodiagonal}, a ``closeness to the diagonal" result is obtained for the optimal alignments. In addition to its own interest, this transversal behavior result is an important technical tool in the proof of our main result. The main Section \ref{sec:mainresult} states the central limit theorem for $L_{n}$ and, after introducing the necessary background on Stein's method, provides its proof. The paper concludes with Section~\ref{sec:Conclusion}, where some pointers to possible future directions of related research are indicated.

\section{Rate of Convergence}\label{sec:RateConvergence}

Our goal below is to prove the following theorem which extends, to multiple sequences and a general scoring function, the work of Alexander~\cite{Alexander:1994}.
\begin{theorem}\label{thm:RateConvOptScoreL1}
Within the setting and with the notation of the introduction,
\begin{align*}
n\gamma^{*}-K\sqrt{n\ln n}\leq\mathbb{E}\left(L_{n}\right)\leq n\gamma^{*},\quad\text{for }\,n\in\mathbb{N},
\end{align*}
where $K$ is a constant independent of $n$.
\end{theorem}

The proof of this first result is based on an adaptation of Rhee's~\cite{Rhee:1995} approach in obtaining Alexander's result. We start by giving the following concentration inequality which immediately follows from Hoeffding's martingale exponential inequality.
\begin{proposition}\label{prop:ConcentLn}
For any $t>0$,
\begin{align}\label{eq:ConcentLn}
\max\left\{\mathbb{P}\left(L_{n}-\mathbb{E}(L_{n})\geq t\right),\,\mathbb{P}\left(L_{n}-\mathbb{E}(L_{n})\leq -t\right)\right\}\leq\exp\left(-\frac{2t^{2}}{nmD^{2}}\right).
\end{align}
\end{proposition}

We now move to the proof of Theorem \ref{thm:RateConvOptScoreL1}. The upper bound is immediate from the super-additivity of $L_{n}$. To obtain the lower bound, we first establish, in the forthcoming proposition, a probabilistic estimate. To do so, for any integer intervals $B_{1},\ldots,B_{m}$, let, below, $L((X^{(1)}_{i})_{i\in B_{1}};(X^{(2)}_{i})_{i\in B_{2}}\,\cdots;(X^{(m)}_{i})_{i\in B_{m}})$ denote the optimal score of $(X^{(j)}_{i})_{i\in B_{j}}$, $j=1,\ldots,m$, and let $[n]:=\{1,\ldots,n\}$, for any $n\in\mathbb{N}$.
\begin{proposition}\label{prop:RheeExten}
For any $m,n\in\mathbb{N}$, and any $x>0$,
\begin{align*}
\mathbb{P}\left(L_{m^{2}n}\geq m^{2}x\right)\leq 2^{-m}m^{4m}n^{2m}\left(\mathbb{P}\left(L_{mn}\geq m\left(x-s^{*}\right)\right)\right)^{1/m}.
\end{align*}
\end{proposition}
\textbf{Proof:} We first claim that, if $L_{m^{2}n}\geq m^{2}x$, we can find $m$ partitions of $[m^{2}n]$:
\begin{align}\label{eq:partition}
\mathbf{B}^{(1)}:=\left(B_{1}^{(1)},B_{2}^{(1)},\ldots,B_{m^{2}}^{(1)}\right),\,\,\ldots\,\,,\,\,\mathbf{B}^{(m)}:=\left(B_{1}^{(m)},B_{2}^{(m)},\ldots,B_{m^{2}}^{(m)}\right),
\end{align}
so that for any $\ell=1,\ldots,m^{2}$,
\begin{align}\label{eq:CellOptScore}
L\left((X^{(1)}_{i})_{i\in B^{(1)}_{\ell}};\,(X^{(2)}_{i})_{i\in B^{(2)}_{\ell}};\,\cdots;(X^{(m)}_{i})_{i\in B^{(m)}_{\ell}}\right)\geq x-s^{*}.
\end{align}
To see this, let us write $x=ks^{*}+\delta$, where $k\in\mathbb{N}_{0}:=\mathbb{N}\cup\{0\}$ and $\delta\in(0,s^{*})$. Set $T_{0}^{(j)}=0$ for every $j\in[m]$. Next, let $(T_{1}^{(j)})_{j\in[m]}$ be the smallest integers such that
\begin{align*}
ks^{*}\leq L\left(X_{1}^{(1)},\ldots,X_{T_{1}^{(1)}}^{(1)};\,\cdots;X_{1}^{(m)},\ldots,X_{T_{1}^{(m)}}^{(m)}\right)<(k+1)s^{*}.
\end{align*}
Such $T_{1}^{(j)}$'s exist since the maximal possible score is $s^{*}$ and since $L_{m^{2}n}\geq m^{2}x$. In general, for $\ell\geq 2$, let $(T_{\ell}^{(j)})_{j\in[m]}$ be the smallest integers with $T_{\ell}^{(j)}>T_{\ell-1}^{(j)}$ for each $j\in[m]$, such that
\begin{align*}
ks^{*}\leq L\left(X_{T_{\ell-1}^{(1)}+1}^{(1)},\ldots,X_{T_{\ell}^{(1)}}^{(1)};\,\ldots;X_{T_{\ell-1}^{(m)}+1}^{(m)},\ldots,X_{T_{\ell}^{(m)}}^{(m)}\right)<(k+1)s^{*}.
\end{align*}
Note that this process ends at a finite time since $m^{2}n<\infty$. Since $L_{m^{2}n}\geq m^{2}x$, the (random) number of cells is at least
\begin{align*}
\frac{m^{2}x}{(k+1)s^{*}}=\frac{m^{2}\left(ks^{*}+\delta\right)}{(k+1)s^{*}}>m^{2}-1.
\end{align*}
By denoting, for $j=1,\ldots,m$,
\begin{align*}
B_{\ell}^{(j)}:=\left(T_{\ell-1}^{(j)}+1,\ldots,T_{\ell}^{(j)}\right),\,\,\,\ell=1,\ldots,m^{2}-1;\quad B_{m^{2}}^{(j)}:=\left(T_{m^{2}-1}^{(j)}+1,\ldots,n\right),
\end{align*}
we have found partitions as in \eqref{eq:partition} which satisfy \eqref{eq:CellOptScore}.

Moreover, since for all $1\leq j\leq m$,
\begin{align*}
\text{Card}\left(B_{1}^{(j)}\right)+\text{Card}\left(B_{2}^{(j)}\right)+\cdots+\text{Card}\left(B_{m^{2}}^{(j)}\right)=m^{2}n,
\end{align*}
there exists $1\leq\ell_{0}\leq m^{2}$ such that
\begin{align*}
\text{Card}\left(B_{\ell_{0}}^{(1)}\right)+\text{Card}\left(B_{\ell_{0}}^{(2)}\right)+\cdots+\text{Card}\left(B_{\ell_{0}}^{(m)}\right)\leq mn.
\end{align*}
Therefore,
\begin{align*}
\left\{L_{m^{2}n}\geq m^{2}x\right\}\subset\bigcup_{\mathscr{S}}\left\{L\left((X_{i}^{(1)})_{i\in B^{(1)}},\ldots,(X_{i}^{(m)})_{i\in B^{(m)}}\right)\geq x-s^{*}\right\},
\end{align*}
where $\mathscr{S}$ is the collection of all $m$ integer intervals $(B^{(1)},\ldots,B^{(m)})$ with $B^{(j)}\subset[m^{2}n]$ for every $j=1,\ldots,m$, such that $\text{Card}(B^{(1)})+\cdots+\text{Card}(B^{(m)})\leq mn$. Notice that
\begin{align*}
\text{Card}(\mathscr{S})\leq\underbrace{\binom{m^{2}n}{2}\cdots\binom{m^{2}n}{2}}_{\text{m many}}=\left(\frac{m^{2}n\left(m^{2}n-1\right)}{2}\right)^{m}\leq 2^{-m}m^{4m}n^{2m}.
\end{align*}
Hence, by Boole's inequality,
\begin{align*}
\mathbb{P}\left(L_{m^{2}n}\geq m^{2}x\right)\leq 2^{-m}m^{4m}n^{2m}\max_{\mathscr{S}}\mathbb{P}\left(L\left((X_{i}^{(1)})_{i\in B^{(1)}},\ldots,(X_{i}^{(m)})_{i\in B^{(m)}}\right)\geq x-s^{*}\right).
\end{align*}

It remains to estimate the last maximum. To do so, first observe that
\begin{align*}
&\max_{\mathscr{S}}\mathbb{P}\left(L\left((X_{i}^{(1)})_{i\in B^{(1)}},\ldots,(X_{i}^{(m)})_{i\in B^{(m)}}\right)\geq x-s^{*}\right)\\
&\quad=\max_{\substack{t_{1},\ldots,t_{m}\in\mathbb{N} \\ \sum t_{i}\leq mn}}\mathbb{P}\left(L\left(X_{1}^{(1)},\ldots,X_{t_{1}}^{(1)};\,\ldots;\,X_{1}^{(m)},\ldots,X_{t_{m}}^{(m)}\right)\geq x-s^{*}\right).
\end{align*}
Next, set
\begin{align*}
\mathbf{Z}_{1}^{(1)}=\left(X_{1}^{(1)},\ldots,X_{t_{1}}^{(1)}\right),\;\mathbf{Z}_{1}^{(2)}=\left(X_{1}^{(2)},\ldots,X_{t_{2}}^{(2)}\right),\;\cdots\;,\mathbf{Z}_{1}^{(m)}=\left(X_{1}^{(m)},\ldots,X_{t_{m}}^{(m)}\right),
\end{align*}
and for $i,j=2,\ldots,m$, let $\mathbf{Z}_{i}^{(j)}$ be an independent copy of $\mathbf{Z}_{1}^{(j)}$. Since the score function $S$ is permutation-invariant,
\begin{align*}
L\left(\mathbf{Z}_{1}^{(1)};\mathbf{Z}_{1}^{(2)};\cdots;\mathbf{Z}_{1}^{(m)}\right)\eqdist L\left(\mathbf{Z}_{2}^{(2)};\mathbf{Z}_{2}^{(3)};\cdots;\mathbf{Z}_{2}^{(1)}\right)\eqdist\cdots\eqdist L\left(\mathbf{Z}_{m}^{(m)};\mathbf{Z}_{m}^{(1)};\cdots;\mathbf{Z}_{m}^{(m-1)}\right).
\end{align*}
Hence,
\begin{align*}
&\left(\mathbb{P}\left(L\left(\mathbf{Z}_{1}^{(1)};\mathbf{Z}_{1}^{(2)};\cdots;\mathbf{Z}_{1}^{(m)}\right)\geq x-s^{*}\right)\right)^{m}\\
&\quad=\prod_{j=1}^{m}\mathbb{P}\left(L\left(\mathbf{Z}_{j}^{(j)};\mathbf{Z}_{j}^{(j+1)};\cdots;\mathbf{Z}_{j}^{(m)};\mathbf{Z}_{j}^{(1)};\cdots;\mathbf{Z}_{j}^{(j-1)}\right)\geq x-s^{*}\right)\\
&\quad\leq\mathbb{P}\left(L\left(\mathbf{Z}_{1}^{(1)}\mathbf{Z}_{2}^{(2)}\ldots\mathbf{Z}_{m}^{(m)};\,\mathbf{Z}_{1}^{(2)}\mathbf{Z}_{2}^{(3)}\ldots\mathbf{Z}_{m}^{(1)};\,\cdots;\,\mathbf{Z}_{1}^{(m)}\mathbf{Z}_{2}^{(1)}\ldots\mathbf{Z}_{m}^{(m-1)}\right)\geq m\left(x-s^{*}\right)\right)\\
&\quad\leq\mathbb{P}\left(L_{mn}\geq m\left(x-s^{*}\right)\right),
\end{align*}
where the last inequality follows by noting that each one of the $m$ sequences in the second to last expression has the same length as $\mathbf{Z}_{1}^{(1)}\mathbf{Z}_{1}^{(2)}\ldots\mathbf{Z}_{1}^{(m)}$, which itself has at most $mn$ terms.\hfill $\Box$

\medskip
The following is essential in establishing the lower bound in Theorem~\ref{thm:RateConvOptScoreL1}.
\begin{proposition}\label{prop:EstExpLn}
For any fixed $m\geq 2$, let $K_{2}\geq m^{2}+2^{-m-3/2}m^{4m+5/2}\sqrt{\pi}\,D$, then
\begin{align*}
\mathbb{E}(L_{m^{2}n})\leq m\,\mathbb{E}(L_{mn})+K_{2}\left(\sqrt{n\ln n}\vee 1\right),\quad\text{for any }\,n\in\mathbb{N}.
\end{align*}
\end{proposition}
\textbf{Proof:} By Proposition \ref{prop:ConcentLn}, for any $z>0$,
\begin{align*}
\mathbb{P}\left(L_{mn}\geq z\right)\leq\exp\left(-\frac{2\left(z-\mathbb{E}(L_{mn})\right)^{2}}{nm^{2}D^{2}}\right).
\end{align*}
Together with Proposition \ref{prop:RheeExten}, for $x>m^{2}s^{*}$,
\begin{align}
\mathbb{P}\left(L_{m^{2}n}\geq x\right)&\leq 2^{-m}m^{4m}n^{2m}\left(\mathbb{P}\left(L_{mn}\geq\frac{x}{m}-ms^{*}\right)\right)^{1/m}\nonumber\\
\label{eq:TailLm2n} &\leq 2^{-m}m^{4m}n^{2m}\exp\left(-\frac{2\left(x-m^{2}s^{*}- m\,\mathbb{E}(L_{mn})\right)^{2}}{nm^{5}D^{2}}\right).
\end{align}

Next, for any $u>0$,
\begin{align}\label{eq:ExpanExpLm2n}
\mathbb{E}(L_{m^{2}n})=\int_{0}^{\infty}\mathbb{P}\left(L_{m^{2}n}\geq x\right)dx\leq u+\int_{u}^{\infty}\mathbb{P}\left(L_{m^{2}n}\geq x\right)dx.
\end{align}
For $n\geq 2$, choosing
\begin{align*}
u:=m\,\mathbb{E}(L_{mn})+m^{2}s^{*}+K_{1}\sqrt{n\ln n},
\end{align*}
where $K_{1}>0$ is a constant, only depending on $m$, to be determined later and using the previous estimate, we obtain that
\begin{align}
\mathbb{E}(L_{m^{2}n})&\leq m\,\mathbb{E}(L_{mn})+m^{2}s^{*}+K_{1}\sqrt{n\ln n}\nonumber\\
\label{eq:ExpLm2n} &\quad\,+\frac{m^{4m}n^{2m}}{2^{m}}\!\int_{m\mathbb{E}(L_{mn})+m^{2}s^{*}+K_{1}\sqrt{n\ln n}}^{\infty}\exp\left(-\frac{2\left(x-m^{2}s^{*}-m\,\mathbb{E}(L_{mn})\right)^{2}}{nm^{5}D^2}\right)dx.
\end{align}
By changing variables, $v=\sqrt{2}(x-m^{2}s^{*}-m\mathbb{E}(L_{mn}))/(\sqrt{n}m^{5/2}D)$, and using
\begin{align*}
\int_{v_{0}}^{\infty}e^{-v^{2}}dv\leq\frac{1}{2v_{0}}\int_{v_{0}}^{\infty}2ve^{-v^{2}}dv=\frac{1}{2v_{0}}e^{-v_{0}^{2}},\quad v_{0}>0,
\end{align*}
the last integral term in \eqref{eq:ExpLm2n} becomes
\begin{align*}
\frac{m^{4m}n^{2m}}{2^{m}}\,\frac{\sqrt{n}m^{5/2}D}{\sqrt{2}}\int_{K_{1}\sqrt{2\ln n}/(m^{5/2}D)}^{\infty}e^{-v^{2}}\,dv\leq\frac{D^2}{2^{m+2}K_{1}\sqrt{\ln n}}\,m^{4m+5}n^{2m+\frac{1}{2}-\frac{2K_{1}^{2}}{m^{5}D^{2}}}.
\end{align*}
Together with \eqref{eq:ExpLm2n}, we have
\begin{align*}
\mathbb{E}\left(L_{m^{2}n}\right)\leq m\,\mathbb{E}(L_{mn})+m^{2}s^{*}+K_{1}\sqrt{n\ln n}+\frac{D^2}{2^{m+2}K_{1}\sqrt{\ln n}}\,m^{4m+5}n^{2m+\frac{1}{2}-\frac{2K_{1}^{2}}{m^{5}D^{2}}}.
\end{align*}
By choosing $K_{1}=m^{3}D$, and letting $K\geq m^{2}s^{*}+m^{3}D+m^{4m+2}D/(2^{m+2})$, it follows that
\begin{align*}
\mathbb{E}(L_{m^{2}n})\leq m\,\mathbb{E}(L_{mn})+K\sqrt{n\ln n},\quad\text{for any }\,n\geq 2.
\end{align*}
When $n=1$, by \eqref{eq:TailLm2n}) and \eqref{eq:ExpanExpLm2n},
\begin{align*}
\mathbb{E}(L_{m^{2}})\leq m\,\mathbb{E}(L_{m})+m^{2}s^{*}+2^{-m-3/2}m^{4m+5/2}D\sqrt{\pi},
\end{align*}
which complete the proof.\hfill $\Box$

\medskip
\noindent
\textbf{Proof of Theorem \ref{thm:RateConvOptScoreL1}:} By Proposition \ref{prop:EstExpLn}, for any $k\in\mathbb{N}$,
\begin{align*}
\frac{\mathbb{E}(L_{m^{k+1}n})}{m^{k+1}n}\leq\frac{\mathbb{E}(L_{m^{k}n})}{m^{k}n}+\frac{K_{2}}{m^{2}}\sqrt{\frac{\ln(m^{k-1}n)}{m^{k-1}n}},\quad\text{for all }\,n\geq 2.
\end{align*}
Summing these inequalities for $1\leq k\leq\ell$ shows that
\begin{align*}
\frac{\mathbb{E}(L_{m^{\ell+1}n})}{m^{\ell+1}n}\leq\frac{\mathbb{E}(L_{mn})}{mn}+\frac{K_{3}}{m}\sqrt{\frac{\ln n}{n}},
\end{align*}
where
\begin{align*}
K_{3}:=\frac{K_{2}}{m}\sum_{k=1}^{\infty}\sqrt{\frac{2(k-1)\ln m+1}{m^{k-1}}}.
\end{align*}
Then, letting $\ell\rightarrow\infty$, leads to
\begin{align}\label{eq:LowBoundExpLmn1}
mn\gamma^{*}\leq\mathbb{E}(L_{mn})+K_{3}\sqrt{n\ln n}\leq\mathbb{E}(L_{mn})+K_{3}\sqrt{mn\ln(mn)},\quad\text{for all }\,n\geq 2.
\end{align}
When $n=1$, Proposition \ref{prop:EstExpLn}, together with \eqref{eq:LowBoundExpLmn1}, implies that
\begin{align}\label{eq:LowBoundExpLmn2}
\mathbb{E}(L_{m})\geq\frac{\mathbb{E}(L_{m^{2}})}{m}-\frac{K_{2}}{m}\geq m\gamma^{*}-K_{3}\sqrt{2\ln m}-\frac{K_{2}}{m}\geq m\gamma^{*}-(K_{2}+K_{3})\sqrt{m\ln m}.
\end{align}
Combining \eqref{eq:LowBoundExpLmn1} and \eqref{eq:LowBoundExpLmn2}, therefore shows that
\begin{align}\label{eq:LowBoundExpLmk}
\mathbb{E}(L_{n})\geq n\gamma^{*}-K\sqrt{n\ln n},\quad\text{for }\,n=mk,\,\,\,k\in\mathbb{N},
\end{align}
where $K=K_{2}+K_{3}$. The end of the proof (when $n$ is not a multiple of $m$) follows from the monotonicity of $\mathbb{E}(L_{n})$.\hfill $\Box$
\begin{remark}\label{rem:Generalizations}
The results obtained to date admit easy generalizations, some of which are described as below and continue to hold throughout the text.
\begin{itemize}
\item [(i)] In \eqref{eq:BddCondScore} (and therefore also in \eqref{eq:bddconditionK}), $D$ can be replaced by $D_{j}$ if $(x_{1},\ldots,x_{m})$ and $(y_{1},\ldots,y_{m})\in\mathcal{A}^{m}$ differ at the $j$-th coordinate, $j=1,\ldots,m$. Accordingly, $mD^{2}$ on the right-hand side of \eqref{eq:ConcentLn} will be replaced by $\sum_{j=1}^{m}D_{j}^{2}$.
\item [(ii)] It is possible to relax the identical distribution assumption between sequences while preserving the independence. More precisely, the convergence of the first moment of $L_{n}$ as well as the rate of convergence remain valid if we only assume that $\mathbf{X}^{(j)}_{n}=(X^{(j)}_{1},\ldots,X^{(j)}_{n})$, $j=1,\ldots,m$, are independent sequences, and that for each $j\in[m]$, $X^{(j)}_{1},\ldots,X^{(j)}_{n}$ have the same distribution, but allow different distributions for random variables in distinct sequences.
\item [(iii)] The alphabet $\mathcal{A}$ can be chosen infinite countable. This, in particular, allows one to study random words with letters taking values in unbounded sets. For example, one could take $m=2$, $\mathcal{A}=\mathbb{N}$ and $S(x,y)=1-(|x-y|/(1+|x-y|))$, to compare random words on unbounded alphabets, and to obtain the same convergence rate for the corresponding first moment.
\end{itemize}
\end{remark}
\begin{remark}\label{rem:NonTrivLimitInfAlphabet}
When $\mathcal{A}$ is countably infinite, one still has $\gamma^{*}>0$ as easily shown by an argument similar to the one provided in the Introduction. Moreover, in the LCS case with a countably infinite alphabet, i.e., $S(x_{1},\cdots,x_{m})={\bf 1}_{\{x_{1}=\cdots=x_{m}\}}$ and $s^{*}=1$, one can show that $\gamma^{*}<1$. Indeed, assume without loss of generality that $\mathcal{A}=\mathbb{N}$, define $m$ new independent sequences $\widetilde{\mathbf{X}}^{(j)}_{n}=(\widetilde{X}^{(j)}_{1},\ldots,\widetilde{X}^{(j)}_{n})$, $j=1,\ldots,m$, by
\begin{align*}
\widetilde{X}_{i}^{(j)}=\left\{\begin{array}{ll} 1,&\text{if }\,X_{i}^{(j)}=1 \\ 2,&\text{if }\,X_{i}^{(j)}\neq 1 \end{array}\right.,\quad i=1,\ldots,n,\quad j=1,\ldots,m.
\end{align*}
Clearly, $L_{n}(\mathbf{X}^{(1)}_{n};\cdots;\mathbf{X}^{(m)}_{n})\leq L_{n}(\widetilde{\mathbf{X}}^{(1)}_{n};\cdots;\widetilde{\mathbf{X}}^{(m)}_{n})$, and thus
\begin{align*}
\lim_{n\rightarrow\infty}\frac{1}{n}\mathbb{E}\left(L_{n}\left(\mathbf{X}^{(1)}_{n};\cdots;\mathbf{X}^{(m)}_{n}\right)\right)\leq\lim_{n\rightarrow\infty}\frac{1}{n}\mathbb{E}\left(L_{n}\left(\widetilde{\mathbf{X}}^{(1)}_{n};\cdots;\widetilde{\mathbf{X}}^{(m)}_{n}\right)\right)<1,
\end{align*}
where the last inequality follows from the fact that the limit is strictly less than 1 for the LCS of multi-sequences with finite alphabet, which can be verified using a multi-sequence analog of~\cite[Lemma 1, Theorem 1]{ChvatalSankoff:1975}.
\end{remark}
\begin{remark}\label{rem:LPP}
Our sequence alignment problem can be represented as a directed last-passage percolation (LPP) problem with dependent weights. More precisely, let the set of vertices be $V:=\{0,1,\ldots,n\}^{m}$, and let the set of oriented edges $\mathcal{E}$ contain all horizontal, vertical and diagonal edges. The horizontal edges are oriented to the right, while the vertical edges are oriented upwards, both have unit lengths. The diagonal edges point up-right at a $45$-degree angle and have length $\sqrt{2}$. With each of the horizontal and vertical edges, we associate a weight of $0$. With the diagonal edge from $(i_{1},\ldots,i_{m})$ to $(i_{1}+1,\ldots,i_{m}+1)$, we associate a weight $S(X_{i_{i}+1}^{(1)},\ldots,X_{i_{m}+1}^{(m)})$ provided that it is strictly positive, and $0$ otherwise. In this manner, we obtain that $L_{n}$ is equal to the total weights of the heaviest path going from $(0,0,\ldots,0)$ to $(n,n,\ldots,n)$ along the lattice.

Another directed LPP representation can be obtained via
\begin{align*}
L_{n}=\max_{\boldsymbol{\pi}\in\text{SI}}\sum_{(i_{1},\ldots,i_{m})\in\boldsymbol{\pi}}S\left(X_{i_{1}}^{(1)},\ldots,X_{i_{m}}^{(m)}\right),
\end{align*}
where SI refers to the set of all \emph{strictly increasing} paths, i.e., paths with all coordinates strictly increasing from a step to another, from $(0,0,\ldots,0)$ to any boundary $x_{j}=n$, $j=1,\ldots,m$.
\end{remark}

\section{Closeness to the Diagonal}\label{sec:closenesstodiagonal}

The purpose of this section is to provide a generalization of~\cite[Theorem 2.1]{HoudreMatzinger:2014} to a multivariate score function setting. For $p_{2},\ldots,p_{m}>0$, let
\begin{align*}
\gamma_{n}(p_{2},\ldots,p_{m}):=\frac{\mathbb{E}\left(L\left(X_{1}^{(1)},\ldots,X_{n}^{(1)};X_{1}^{(2)},\ldots,X_{np_{2}}^{(2)};\cdots;X_{1}^{(m)},\ldots,X_{np_{m}}^{(m)}\right)\right)}{n\left(1+\sum_{i=2}^{m}p_{i}\right)/m},
\end{align*}
where, when not integers, the indices $np_{j}$, $j=2,\ldots,m$, are understood to be rounded-up to the nearest positive integers, and let
\begin{align*}
\gamma^{*}(p_{2},\ldots,p_{m}):=\lim_{n\rightarrow\infty}\gamma_{n}(p_{2},\ldots,p_{m}).
\end{align*}
Moreover, for $q_{1},\ldots,q_{m}>0$ with $q_{1}+\cdots+q_{m}=m$, let
\begin{align*}
\widetilde{\gamma}_{n}(q_{1},\ldots,q_{m}):=\frac{1}{n}\mathbb{E}\left(L\left(X_{1}^{(1)},\ldots,X_{nq_{1}}^{(1)};\cdots;X_{1}^{(m)},\ldots,X_{nq_{m}}^{(m)}\right)\right),
\end{align*}
and let
\begin{align*}
\widetilde{\gamma}^{*}(q_{1},\ldots,q_{m}):=\lim_{n\rightarrow\infty}\widetilde{\gamma}_{n}(q_{1},\ldots,q_{m}).
\end{align*}
The function $\gamma^{*}$ is a re-parametrization of $\widetilde{\gamma}^{*}$, i.e.,
\begin{align*}
\gamma^{*}(p_{2},\ldots,p_{m})=\widetilde{\gamma}^{*}\left(q_{1}(p_{2},\ldots,p_{m}),\ldots,q_{m}(p_{2},\ldots,p_{m})\right),
\end{align*}
with
\begin{align}\label{eq:Transformpjqj}
q_{j}(p_{2},\ldots,p_{m})=\frac{mp_{j}}{1+\sum_{k=2}^{m}p_{k}},\quad j=2,\ldots,m.
\end{align}
A super-additivity argument, as in~\cite{ChvatalSankoff:1975}, shows that both limits above do exist, and depend on the size of the alphabet, the score function $S$, as well as the distribution of $X_{1}^{(1)}$, but none of this is of importance for our purposes.
\begin{proposition}\label{prop:ConcavTildeGamma}
The function $\widetilde{\gamma}^{*}$ is concave on its domain. In particular, $\widetilde{\gamma}^{*}$ achieves its maximum at $q_{1}=\cdots=q_{m}=1$.
\end{proposition}
\textbf{Proof:} For any $c_{1},c_{2}>0$ with $c_{1}+c_{2}=1$, any $q_{j},r_{j}>0$, $j=1,\ldots,m$, with $\sum_{j=1}^{m}q_{j}=m$ and $\sum_{j=1}^{m}r_{j}=m$, we have
\begin{align*}
\widetilde{\gamma}_{n}\left(c_{1}(q_{1},\ldots,q_{m})+c_{2}(r_{1},\ldots,r_{m})\right)&=\widetilde{\gamma}_{n}\left(c_{1}q_{1}+c_{2}r_{1},\ldots,c_{1}q_{m}+c_{2}r_{m}\right)\\
&=\frac{1}{n}\mathbb{E}\left(L\left(\left(X_{i}^{(1)}\right)_{i=1}^{n(c_{1}q_{1}+c_{2}r_{1})};\cdots;\left(X_{i}^{(m)}\right)_{i=1}^{n(c_{1}q_{m}+c_{2}r_{m})}\right)\right)\\
&\geq c_{1}\,\frac{1}{c_{1}n}\mathbb{E}\left(L\left(X_{1}^{(1)},\ldots,X_{c_{1}nq_{1}}^{(1)};\cdots;X_{1}^{(m)},\ldots,X_{c_{1}nq_{m}}^{(m)}\right)\right)\\
&\quad\,\!+c_{2}\,\frac{1}{c_{2}n}\mathbb{E}\!\left(L\!\left(X_{1}^{(1)},\ldots,X_{c_{2}nr_{1}}^{(1)};\cdots;X_{1}^{(m)},\ldots,X_{c_{2}nr_{m}}^{(m)}\right)\right)\\
&=c_{1}\widetilde{\gamma}_{c_{1}n}(q_{1},\ldots,q_{m})+c_{2}\widetilde{\gamma}_{c_{2}n}(r_{1},\ldots,r_{m}).
\end{align*}
The concavity of $\widetilde{\gamma}^{*}$ follows immediately by taking limits, as $n\rightarrow\infty$, on both sides of the above inequality.

Moreover, since the score function $S$ is permutation-invariant so is $\widetilde{\gamma}^{*}$. Now assume that there exist $q_{1},\ldots,q_{m}>0$ with $\sum_{j=1}^{m}q_{j}=m$, such that $\widetilde{\gamma}^{*}(q_{1},\ldots,q_{m})>\widetilde{\gamma}^{*}(1,\ldots,1)$. Let $(q_{1}^{(\ell)},\ldots,q_{m}^{(\ell)})$, $\ell=1,\ldots,m!$, be all the permutations of $(q_{1},\ldots,q_{m})$, then
\begin{align*}
\sum_{\ell=1}^{m!}q_{k}^{(\ell)}=m!,\quad\text{for all }\,k=1,\ldots,m,
\end{align*}
and
\begin{align*}
\widetilde{\gamma}^{*}\left(q_{1}^{(1)},\ldots,q_{m}^{(1)}\right)=\cdots=\widetilde{\gamma}^{*}\left(q_{1}^{(m!)},\ldots,q_{m}^{(m!)}\right)>\widetilde{\gamma}^{*}(1,\ldots,1).
\end{align*}
It then follows that
\begin{align*}
\widetilde{\gamma}^{*}\left(\frac{1}{m!}\sum_{\ell=1}^{m!}\left(q_{1}^{(\ell)},\ldots,q_{m}^{(\ell)}\right)\right)=\widetilde{\gamma}^{*}(1,\ldots,1)<\widetilde{\gamma}^{*}(q_{1},\ldots,q_{m})=\frac{1}{m!}\sum_{\ell=1}^{m!}\widetilde{\gamma}^{*}\left(q_{1}^{(\ell)},\ldots,q_{m}^{(\ell)}\right),
\end{align*}
which contradicts the concavity property of $\widetilde{\gamma}^{*}$.\hfill $\Box$

\medskip
Therefore, $\gamma^{*}(p_{2},\ldots,p_{m})$ reaches its maximum at $p_{2}=\cdots=p_{m}=1$, which is equal to $\gamma^{*}$ and is given as in \eqref{eq:LimitMeanLn}. Moreover, since $\widetilde{\gamma}^{*}$ is concave, there exists $\delta^{*}\in(0,\gamma^{*})$, such that
\begin{align}
\{(1,\ldots,1)\}&\varsubsetneq\left\{(q_{2},\ldots,q_{m}):\,\, q_{j}>0,\,j=2,\ldots,m;\,\,\widetilde{\gamma}^{*}\left(1-\sum_{j=2}^{m}q_{j},q_{2},\ldots,q_{m}\right)\ge\delta^{*}\right\}\nonumber\\
\label{eq:LevelSetdelta} &\varsubsetneq\left\{ (q_{2},\ldots,q_{m}):\,\,q_{j}>0,\,j=2,\ldots,m;\,\,\sum_{j=2}^{m}q_{j}<m\right\}.
\end{align}
That is, the upper level set of $\widetilde{\gamma}^{*}$ at level $\delta$ is a pure subset of its domain, and contains other points than the singleton $\{(1,\ldots,1)\}$. Under the transformation \eqref{eq:Transformpjqj}, the inverse image of the above level set is a closed and bounded subset of $\{(p_{2},\ldots,p_{m}):p_{2},\ldots,p_{m}>0\}$, which contains $(1,\ldots,1)$. We can then choose $\boldsymbol{p}^{(1)}=(p_{1}^{(1)},\ldots,p_{m}^{(1)})$ and $\boldsymbol{p}^{(2)}=(p_{1}^{(2)},\ldots,p_{m}^{(2)})$, with $0<p_{j}^{(1)}<1<p_{j}^{(2)}$, $j=2,\ldots,m$, so that
\begin{align*}
\left[p_{2}^{(1)},\,p_{2}^{(2)}\right]\times\cdots\times\left[p_{m}^{(1)},\,p_{m}^{(2)}\right]
\end{align*}
is the smallest rectangle containing this inverse image. Hence, for any $\delta\in[\delta^{*},\gamma^{*})$, and any $p_{2},\ldots,p_{m}>0$ such that $p_{j}<p_{j}^{(1)}$ or $p_{j}>p_{j}^{(2)}$ for at least one $j\in\{2,\ldots,m\}$, we have
\begin{align}\label{eq:GammaDelta}
\gamma^{*}(p_{2},\ldots,p_{m})\leq\delta<\gamma^{*}.
\end{align}

Let us now introduce some further notations. Assume that $n=vd$ for some $v,d\in\mathbb{N}$, and for every $j=2,\ldots,m$, let the integers
\begin{align}\label{eq:OptCond1}
0=r_{0}^{(j)}\leq r_{1}^{(j)}\leq\cdots\leq r_{d-1}^{(j)}\leq r_{d}^{(j)}=n,
\end{align}
be such that
\begin{align}\label{eq:OptCond2}
L_{n}\!\left(\mathbf{X}_{n}^{(1)};\,\cdots;\mathbf{X}_{n}^{(m)}\right)\!=\!\sum_{k=1}^{d}L\!\left(\!\!\left(X_{i}^{(1)}\right)_{i=v(k-1)+1}^{vk};\left(X_{i}^{(2)}\right)_{i=r_{k-1}^{(2)}+1}^{r_{k}^{(2)}};\cdots;\left(X_{i}^{(m)}\right)_{i=r_{k-1}^{(m)}+1}^{r_{k}^{(m)}}\!\right).
\end{align}
For $\varepsilon>0$, let $E_{n,\varepsilon}(\boldsymbol{p}^{(1)},\boldsymbol{p}^{(2)})$ be the event that for all $(m-1)$-tuples of integer vectors
$(r_{1}^{(j)},\ldots,r_{d-1}^{(j)})$, $j=2,\ldots,m$, satisfying \eqref{eq:OptCond1} and \eqref{eq:OptCond2}, we have
\begin{align}\label{eq:OptCond3}
\text{Card}\left(\left\{k\in\{1,\ldots,d\}:\,vp_{j}^{(1)}\leq r_{k}^{(j)}-r_{k-1}^{(j)}\leq vp_{j}^{(2)},\,\,\forall j=2,\ldots,m\right\}\right)\geq (1-\varepsilon)\,d.
\end{align}
Also, let $\mathcal{R}_{n,\varepsilon}(\boldsymbol{p}^{(1)},\boldsymbol{p}^{(2)})$ be the set of all non-random $(m-1)$-tuples of integer vectors $(r_{1}^{(j)},\ldots,r_{d-1}^{(j)})$, $j=2,\ldots,m$, satisfying \eqref{eq:OptCond1} and \eqref{eq:OptCond3}, and let
$\overline{\mathcal{R}}_{n,\varepsilon}(\boldsymbol{p}^{(1)},\boldsymbol{p}^{(2)})$ be the set of all non-random $(m-1)$-tuples of integer vectors
$(r_{1}^{(j)},\ldots,r_{d-1}^{(j)})$, $j=2,\ldots,m$, satisfying \eqref{eq:OptCond1}, but not \eqref{eq:OptCond3}. We begin with the following lemma.
\begin{lemma}\label{lem:EstiExpDiff}
Let $\varepsilon>0$, let $\delta^{*}>0$ and let the vectors $\boldsymbol{p}^{(1)},\boldsymbol{p}^{(2)}$ be chosen, as above, so that \eqref{eq:LevelSetdelta} and \eqref{eq:GammaDelta} are satisfied. For any vectors $(r_{1}^{(j)},\ldots,r_{d-1}^{(j)})$, $j=2,\ldots,m$, in $\overline{\mathcal{R}}_{n,\varepsilon}(\boldsymbol{p}^{(1)},\boldsymbol{p}^{(2)})$, and any $\eta\in(0,\gamma^{*}-\delta^{*})$,
\begin{align*}
\mathbb{E}\left(\sum_{k=1}^{d}L\left(\left(X_{i}^{(1)}\right)_{i=v(k-1)+1}^{vk};\left(X_{i}^{(2)}\right)_{i=r_{k-1}^{(2)}+1}^{r_{k}^{(2)}};\cdots;\left(X_{i}^{(m)}\right)_{i=r_{k-1}^{(m)}+1}^{r_{k}^{(m)}}\!\right)-L_{n}\right)\leq -\frac{\eta\,\varepsilon n}{m}.
\end{align*}
for all $n=n(\varepsilon,\eta,\boldsymbol{p}^{(1)},\boldsymbol{p}^{(2)})$ large enough.
\end{lemma}
\textbf{Proof:} By super-additivity, for any $p_{2},\ldots,p_{m}>0$ and any $n\in\mathbb{N}$,
\begin{align*}
\gamma^{*}(p_{2},\ldots,p_{m})&=\lim_{n\rightarrow\infty}\frac{\mathbb{E}\left(L\left(X_{1}^{(1)},\ldots,X_{n}^{(1)};X_{1}^{(2)},\ldots,X_{np_{2}}^{(2)};\cdots;X_{1}^{(m)},\ldots,X_{np_{m}}^{(m)}\right)\right)}{n\left(1+\sum_{i=2}^{m}p_{i}\right)/m}\\
&\geq\frac{m\,\mathbb{E}\left(L\left(X_{1}^{(1)},\ldots,X_{n}^{(1)};X_{1}^{(2)},\ldots,X_{np_{2}}^{(2)};\cdots;X_{1}^{(m)},\ldots,X_{np_{m}}^{(m)}\right)\right)}{n\left(1+\sum_{i=2}^{m}p_{i}\right)}
\end{align*}
Together with \eqref{eq:GammaDelta}, for any $p_{2},\ldots,p_{m}>0$ such that $p_{j}<p_{j}^{(1)}$ or $p_{j}>p_{j}^{(2)}$ for at least one $j\in\{2,\ldots,m\}$,
\begin{align*}
\frac{m\,\mathbb{E}\left(L\left(X_{1}^{(1)},\ldots,X_{v}^{(1)};X_{1}^{(2)},\ldots,X_{vp_{2}}^{(2)};\cdots;X_{1}^{(m)},\ldots,X_{vp_{m}}^{(m)}\right)\right)}{v\left(1+\sum_{i=2}^{m}p_{i}\right)}\leq\delta^{*}.
\end{align*}
Since $\mathbf{X}_{n}^{(j)}$, $j=1,\ldots,m$, are sequences of i.i.d. random variables, for each $k=1,\ldots,d$, by taking $vp_{j}=r_{k}^{(j)}-r_{k-1}^{(j)}$, if $(r_{k}^{(j)}-r_{k-1}^{(j)})/v\notin[p_{j}^{(1)},p_{j}^{(2)}]$ for at least one $j\in\{2,\ldots,m\}$, the above inequality becomes
\begin{align*}
\gamma^{*}-\frac{m\,\mathbb{E}\left(L\left(\left(X_{i}^{(1)}\right)_{i=v(k-1)+1}^{vk};\left(X_{i}^{(2)}\right)_{i=r_{k-1}^{(2)}+1}^{r_{k}^{(2)}};\cdots;\left(X_{i}^{(m)}\right)_{i=r_{k-1}^{(m)}+1}^{r_{k}^{(m)}}\right)\right)}{v+\sum_{i=2}^{m}\left(r_{k}^{(j)}-r_{k-1}^{(j)}\right)}\geq\gamma^{*}-\delta^{*},
\end{align*}
which implies that
\begin{align*}
&\frac{\gamma^{*}}{m}\!\left[v\!+\!\!\sum_{i=2}^{m}\!\left(r_{k}^{(j)}\!-\!r_{k-1}^{(j)}\right)\!\right]\!-\!\mathbb{E}\!\left(\!L\!\left(\!\left(X_{i}^{(1)}\right)_{i=v(k-1)+1}^{vk}\!;\left(X_{i}^{(2)}\right)_{i=r_{k-1}^{(2)}+1}^{r_{k}^{(2)}}\!;\cdots;\left(X_{i}^{(m)}\right)_{i=r_{k-1}^{(m)}+1}^{r_{k}^{(m)}}\right)\!\right)\\
&\quad\geq\frac{\gamma^{*}-\delta^{*}}{m}\left[v+\sum_{i=2}^{m}\left(r_{k}^{(j)}-r_{k-1}^{(j)}\right)\right]\geq\frac{v(\gamma^{*}-\delta^{*})}{m}.
\end{align*}
Let $\mathcal{M}$ be the collection of indices $k\in\{1,\ldots,d\}$ such that $r_{k}^{(j)}-r_{k-1}^{(j)}\notin[vp_{k}^{(1)},vp_{k}^{(2)}]$ for at
least one $j\in\{2,\ldots,m\}$, then if $(r_{1}^{(j)},\ldots,r_{d-1}^{(j)})$, $j=2,\ldots,m$, belong to $\overline{\mathcal{R}}_{n,\varepsilon}(\boldsymbol{p}^{(1)},\boldsymbol{p}^{(2)})$,
\begin{align}
&\sum_{k\in\mathcal{M}}\!\left[\frac{\gamma^{*}}{m}\!\!\left(\!v\!+\!\!\sum_{i=2}^{m}\!\left(\!r_{k}^{(j)}\!-\!r_{k-1}^{(j)}\!\right)\!\!\right)\!\!-\!\mathbb{E}\!\left(\!L\!\left(\!\!\left(\!X_{i}^{(1)}\!\right)_{i=v(k-1)+1}^{vk}\!;\!\left(\!X_{i}^{(2)}\!\right)_{i=r_{k-1}^{(2)}+1}^{r_{k}^{(2)}}\!;\cdots;\!\left(\!X_{i}^{(m)}\!\right)_{i=r_{k-1}^{(m)}+1}^{r_{k}^{(m)}}\right)\!\!\right)\!\right]\nonumber\\
\label{eq:LowBoundDiff}
&\quad\geq\sum_{k\in\mathcal{M}}\frac{v(\gamma^{*}-\delta^{*})}{m}\geq\frac{v(\gamma^{*}-\delta^{*})}{m}\varepsilon d=\frac{(\gamma^{*}-\delta^{*})\varepsilon n}{m}.
\end{align}
On the other hand,
\begin{align}
&\sum_{k\in\mathcal{M}}\!\left[\frac{\gamma^{*}}{m}\!\!\left(\!v\!+\!\!\sum_{i=2}^{m}\!\left(\!r_{k}^{(j)}\!-\!r_{k-1}^{(j)}\!\right)\!\!\right)\!\!-\!\mathbb{E}\!\left(\!\!L\!\left(\!\!\left(\!X_{i}^{(1)}\!\right)_{i=v(k-1)+1}^{vk}\!;\!\left(\!X_{i}^{(2)}\!\right)_{i=r_{k-1}^{(2)}+1}^{r_{k}^{(2)}}\!;\cdots;\!\left(\!X_{i}^{(m)}\!\right)_{i=r_{k-1}^{(m)}+1}^{r_{k}^{(m)}}\right)\!\!\right)\!\right]\nonumber\\
&\quad\leq\sum_{k=1}^{d}\!\left[\frac{\gamma^{*}}{m}\!\!\left(\!v\!+\!\!\sum_{i=2}^{m}\!\left(\!r_{k}^{(j)}\!-\!r_{k-1}^{(j)}\!\right)\!\!\right)\!\!-\!\mathbb{E}\!\left(\!\!L\!\left(\!\!\left(\!X_{i}^{(1)}\!\right)_{i=v(k-1)+1}^{vk}\!;\left(\!X_{i}^{(2)}\!\right)_{i=r_{k-1}^{(2)}+1}^{r_{k}^{(2)}}\!;\cdots;\left(\!X_{i}^{(m)}\!\right)_{i=r_{k-1}^{(m)}+1}^{r_{k}^{(m)}}\!\right)\!\!\right)\!\right]\nonumber\\
\label{eq:UpBoundDiff} &\quad=n\gamma^{*}-\mathbb{E}\left(\sum_{k=1}^{d}L\!\left(\!\!\left(X_{i}^{(1)}\right)_{i=v(k-1)+1}^{vk};\left(X_{i}^{(2)}\right)_{i=r_{k-1}^{(2)}+1}^{r_{k}^{(2)}};\cdots;\left(X_{i}^{(m)}\right)_{i=r_{k-1}^{(m)}+1}^{r_{k}^{(m)}}\!\right)\right).
\end{align}
Combining \eqref{eq:LowBoundDiff} and \eqref{eq:UpBoundDiff}, gives
\begin{align}\label{eq:LowBoundGammaDiff}
n\gamma^{*}\!-\!\mathbb{E}\!\left(\sum_{k=1}^{d}\!L\!\left(\!\!\left(X_{i}^{(1)}\right)_{i=v(k-1)+1}^{vk}\!;\left(X_{i}^{(2)}\right)_{i=r_{k-1}^{(2)}+1}^{r_{k}^{(2)}}\!;\cdots;\left(X_{i}^{(m)}\right)_{i=r_{k-1}^{(m)}+1}^{r_{k}^{(m)}}\!\right)\!\right)\!\ge\!\frac{(\gamma^{*}\!-\!\delta^{*})\varepsilon n}{m},
\end{align}
as long as $(r_{1}^{(j)},\ldots,r_{d-1}^{(j)})$, $j=2,\ldots,m$, belong to $\overline{\mathcal{R}}_{n,\varepsilon}(\boldsymbol{p}^{(1)},\boldsymbol{p}^{(2)})$.

Now by Theorem \ref{thm:RateConvOptScoreL1}, for $n$ large enough (depending only on $\delta$, $\varepsilon$, $\boldsymbol{p}^{(1)}$ and $\boldsymbol{p}^{(2)}$),
\begin{align}\label{eq:ConvRateDelta}
0\leq\gamma^{*}-\frac{\mathbb{E}(L_{n})}{n}\leq\frac{(\gamma^{*}-\delta^{*}-\eta)\varepsilon}{m}.
\end{align}
Combining \eqref{eq:LowBoundGammaDiff} and \eqref{eq:ConvRateDelta} finishes the proof.\hfill $\Box$

\medskip
The main result of this section, as stated next, shows that the event $E_{n,\varepsilon}(\boldsymbol{p}^{(1)},\boldsymbol{p}^{(2)})$ holds with high probability.
\begin{theorem}\label{thm:CloseDiagonal1}
Let $\varepsilon>0$. Let $\boldsymbol{p}^{(1)}=(p_{1}^{(1)},\ldots,p_{m}^{(1)})$ and $\boldsymbol{p}^{(2)}=(p_{1}^{(2)},\ldots,p_{m}^{(2)})$, with $0<p_{j}^{(1)}<1<p_{j}^{(2)}$, $j=2,\ldots,m$, such that \eqref{eq:GammaDelta} holds. For any $\eta\in(0,\gamma^{*}-\delta^{*})$, fix the integer $v$ to be such that
\begin{align*}
\frac{1+\ln(1+v)}{v}\leq\frac{\eta^{2}\varepsilon^{2}}{2m^{3}(m-1)D^{2}}.
\end{align*}
Then for $n=n(\eta,\varepsilon,\boldsymbol{p}^{(1)},\boldsymbol{p}^{(2)})$ large enough,
\begin{align*}
\mathbb{P}\left(E_{n,\varepsilon}(\boldsymbol{p}^{(1)},\boldsymbol{p}^{(2)})\right)\geq 1-\exp\left(-n\left(-\frac{1+\ln(1+v)}{v}(m-1)+\frac{\eta^{2}\varepsilon^{2}}{2m^{3}D^{2}}\right)\right).
\end{align*}
\end{theorem}
\textbf{Proof:} Clearly,
\begin{align}\label{eq:EstCardbarReps}
\text{Card}\!\left(\overline{\mathcal{R}}_{n,\varepsilon}(\boldsymbol{p}^{(1)},\boldsymbol{p}^{(2)})\right)\!\leq\!\binom{n\!+\!d}{d}^{m-1}\!\!\!\!\!\leq\!\frac{(n\!+\!d)^{d(m-1)}}{(d\,!)^{m-1}}\!\leq\!\left(\!\frac{e(n\!+\!d)}{d}\!\right)^{d(m-1)}\!\!\!\!\!\!\!=\!\left(e(1\!+\!v)\right)^{d(m-1)}.
\end{align}
Fix any $\eta\in(0,\gamma^{*}-\delta^{*})$. For the event $E_{n,\varepsilon}(\boldsymbol{p}^{(1)},\boldsymbol{p}^{(2)})$ not to hold, there must exist at least one element $\vec{r}:=\{(r_{1}^{(j)},\ldots,r_{d-1}^{(j)}),\,j=2,\ldots,m\}$ of $\overline{\mathcal{R}}_{n,\varepsilon}(\boldsymbol{p}^{(1)},\boldsymbol{p}^{(2)})$, defining an optimal alignment, i.e., one element $\vec{r}$ must be such that
\begin{align*}
L(\vec{r}):=\sum_{k=1}^{d}L\left(\left(X_{i}^{(1)}\right)_{i=v(k-1)+1}^{vk};\left(X_{i}^{(2)}\right)_{i=r_{k-1}^{(2)}+1}^{r_{k}^{(2)}};\cdots;\left(X_{i}^{(m)}\right)_{i=r_{k-1}^{(m)}+1}^{r_{k}^{(m)}}\right)\geq L_{n}.
\end{align*}
Hence,
\begin{align*}
\mathbb{P}\left(E_{n,\varepsilon}^{c}(\boldsymbol{p}^{(1)},\boldsymbol{p}^{(2)})\right)\leq\sum_{\vec{r}\in\overline{\mathcal{R}}_{n,\varepsilon}(\boldsymbol{p}^{(1)},\boldsymbol{p}^{(2)})}\mathbb{P}\left(L(\vec{r})-L_{n}\geq 0\right).
\end{align*}
When $\vec{r}\in\overline{\mathcal{R}}_{n,\varepsilon}(\boldsymbol{p}^{(1)},\boldsymbol{p}^{(2)})$, it follows from Lemma \ref{lem:EstiExpDiff} and
Hoeffding's exponential martingale inequality that, for $n$ large enough,
\begin{align*}
\mathbb{P}\left(L(\vec{r})-L_{n}\geq 0\right)\leq\mathbb{P}\left(L(\vec{r})-L_{n}-\mathbb{E}\left(L(\vec{r})-L_{n}\right)\geq\frac{\eta\varepsilon n}{m}\right)\leq\exp\left(-\frac{\eta^{2}\varepsilon^{2}}{2m^{3}D^{2}}n\right).
\end{align*}
This last inequality, together with \eqref{eq:EstCardbarReps}, finally leads to
\begin{align}
\mathbb{P}\left(E_{n,\varepsilon}^{c}(\boldsymbol{p}^{(1)},\boldsymbol{p}^{(2)})\right)&\leq\left(e(1+v)\right)^{d(m-1)}\exp\left(-\frac{\eta^{2}\varepsilon^{2}}{2m^{3}D^{2}}n\right)\nonumber\\
\label{eq:EstProbCompEn} &=\exp\left(-n\left(-\frac{1+\ln(1+v)}{v}(m-1)+\frac{\eta^{2}\varepsilon^{2}}{2m^{3}D^{2}}\right)\right),
\end{align}
which completes the proof.\hfill $\Box$

\begin{remark}\label{rmk:epsilonchoice}
Theorem \ref{thm:CloseDiagonal1} has been shown to hold for sufficiently large $n$. but more quantitative estimates can be obtained using Theorem \ref{thm:RateConvOptScoreL1}. To do so, first recall that by Theorem \ref{thm:RateConvOptScoreL1}, there exists a constant $K>0$, independent of $n$, such that
\begin{equation}\label{eq:scaledgamma}
\gamma^{*}-\frac{\mathbb{E}(L_{n})}{n}\leq K\sqrt{\frac{\ln n}{n}},\quad\text{for all }\,n\geq 1.
\end{equation}
Next, note that the terms related to $n$ in Theorem \ref{thm:CloseDiagonal1} should satisfy
\begin{align*}
\frac{1+\ln(1+v)}{v}\leq\frac{\eta^{2}\varepsilon^{2}}{2m^{3}(m-1)D^{2}}\qquad\text{and}\qquad\gamma^{*}-\frac{\mathbb{E}(L_{n})}{n}\leq\frac{(\gamma^{*}-\delta^{*}-\eta)\varepsilon}{m},
\end{align*}
where $v$ will be taken to be $n^{\alpha}$, $\alpha\in(0,1)$. By \eqref{eq:scaledgamma}, the latter condition will be satisfied if
\begin{align*}
K^{2}m^{2}\ln n\leq\varepsilon^{2}\left(\gamma^{*}-\delta^{*}-\eta\right)n.
\end{align*}

Next, choose $c_{1}\in\mathbb{R}$ so that
\begin{align*}
c_{1}^{2}\geq\frac{4m^{3}(m-1)D^{2}}{\eta^{2}},
\end{align*}
and that
\begin{align*}
c_{1}^{2}\left(\frac{1+\ln\left(n^{\alpha}+1\right)}{n^{\alpha}}\right)\geq\frac{K^{2}m^{2}\ln n}{\left(\gamma^{*}-\delta^{*}-\eta\right)n},\quad\text{for all }\,n\in\mathbb{N}.
\end{align*}
Hence, by taking
\begin{align*}
\varepsilon^{2}:=c_{1}^{2}\left(\frac{1+\ln\left(n^{\alpha}+1\right)}{n^{\alpha}}\right),
\end{align*}
Theorem \ref{thm:CloseDiagonal1} with $v=n^{\alpha}$ provides the estimate
\begin{align}\label{eq:Zequalszeroestimate}
\mathbb{P}\left(E_{n,\varepsilon}(\boldsymbol{p}^{(1)},\boldsymbol{p}^{(2)})\right)\geq 1-e^{-n^{1-\alpha}(m-1)(1+\ln(1+n^{\alpha}))},\quad\text{for all }\,n\in\mathbb{N}.
\end{align}
\end{remark}

Below we will also need another interpretation of Theorem \ref{thm:CloseDiagonal1}, which we discuss next. Let $D_{n,\varepsilon}(\boldsymbol{p}^{(1)},\boldsymbol{p}^{(2)})$ be the event that all the points representing any optimal alignment of
$\mathbf{X}^{(j)}_{n}=(X^{(j)}_{1},\ldots,X^{(j)}_{n})$, $j=1,\ldots,m$, are within the parallelepiped in $\mathbb{R}^{m}$ bounded by
\begin{align*}
x_{j}=p_{j}^{(1)}x_{1}-p_{j}^{(1)}n\varepsilon-p_{j}^{(1)}v,\quad x_{j}=p_{j}^{(2)}x_{1}+p_{j}^{(2)}n\varepsilon+p_{j}^{(2)}v,\quad j=2,\ldots,m.
\end{align*}
The following result measures the ``closeness to the diagonal".
\begin{theorem}\label{thm:CloseDiagonal2}
Let $\varepsilon>0$. Let $\boldsymbol{p}^{(1)}=(p_{1}^{(1)},\ldots,p_{m}^{(1)})$ and $\boldsymbol{p}^{(2)}=(p_{1}^{(2)},\ldots,p_{m}^{(2)})$, with $0<p_{j}^{(1)}<1<p_{j}^{(2)}$, $j=2,\ldots,m$, such that \eqref{eq:GammaDelta} holds. For any $\eta\in(0,\gamma^{*}-\delta^{*})$, fix the integer $v$ to be such that
\begin{align*}
\frac{1+\ln(1+v)}{v}\leq\frac{\eta^{2}\varepsilon^{2}}{2m^{3}(m-1)D^{2}}.
\end{align*}
Then for $n=n(\eta,\varepsilon,\boldsymbol{p}^{(1)},\boldsymbol{p}^{(2)})$ large enough,
\begin{align*}
\mathbb{P}\left(D_{n,\varepsilon}(\boldsymbol{p}^{(1)},\boldsymbol{p}^{(2)})\right)\geq 1-2(m-1)\exp\left(-n\left(-\frac{1+\ln(1+v)}{v}(m-1)+\frac{\eta^{2}\varepsilon^{2}}{2m^{3}D^{2}}\right)\right).
\end{align*}
\end{theorem}
\textbf{Proof:} For $j=2,\ldots,m$, let $D_{1,j}$ be the event that any optimal alignment of $\mathbf{X}^{(1)}_{n},\ldots,\mathbf{X}^{(m)}_{n}$, is ``above" the hyperplane $x_{j}=p_{j}^{(1)}x_{1}-p_{j}^{(1)}n\varepsilon-p_{j}^{(1)}v$, i.e., if the $x_{1}$-th element of $\mathbf{X}^{(1)}_{n}$ is aligned with $x_{j}$-th element of $\mathbf{X}^{(j)}_{n}$, then $x_{j}\geq p_{j}^{(1)}x_{1}-p_{j}^{(1)}n\varepsilon-p_{j}^{(1)}v$. Similarly, for $j=2,\ldots,m$, let $D_{2,j}$ be the event that any optimal alignment of $\mathbf{X}^{(1)}_{n},\ldots,\mathbf{X}^{(m)}_{n}$, is ``below" the
hyperplane $x_{j}=p_{j}^{(2)}x_{1}+p_{j}^{(2)}n\varepsilon+p_{j}^{(2)}v$, i.e., if the $x_{1}$-th element of $\mathbf{X}^{(1)}_{n}$ is aligned with $x_{j}$-th element of $\mathbf{X}^{(j)}_{n}$, then $x_{j}\leq p_{j}^{(2)}x_{1}+p_{j}^{(2)}n\varepsilon+p_{j}^{(2)}v$. Clearly,
\begin{align*}
D_{n,\varepsilon}(\boldsymbol{p}^{(1)},\boldsymbol{p}^{(2)})=\bigcap_{j=2}^{m}\left(D_{1,j}\cap D_{2,j}\right),
\end{align*}
Theorem \ref{thm:CloseDiagonal2} will then be a consequence of Theorem \ref{thm:CloseDiagonal1} together with the following inclusions:
\begin{align}\label{eq:InclEnDn}
E_{n,\varepsilon}(\boldsymbol{p}^{(1)},\boldsymbol{p}^{(2)})\subset D_{1,j},\quad E_{n,\varepsilon}(\boldsymbol{p}^{(1)},\boldsymbol{p}^{(2)})\subset D_{2,j},\quad\text{for all }\,j=2,\ldots,m.
\end{align}

We now verify the first inclusion in \eqref{eq:InclEnDn}, while the second inclusion can be shown in a similar way. Let us first assume that $x_{1}$ is a multiple of $v$: $x_{1}=\ell v$ for some $\ell\in\mathbb{N}$, and that $x_{1}\le n\varepsilon$. Then, for any $j=2,\ldots,m$, $p_{j}^{(1)}x_{1}-p_{j}^{(1)}n\varepsilon\leq 0$. If any alignment (in particular, any optimal alignment) aligns the $x_{1}$-th element of $\mathbf{X}^{(1)}_{n}$ with the $x_{j}$ element of $\mathbf{X}^{(j)}_{n}$, then $x_{j}$ must be an integer between $1$ and $n$. In particular, $x_{j}\ge 0\ge p_{j}^{(1)}x_{1}-p_{j}^{(1)}n\varepsilon$. Next, assume that $x_{1}=\ell v$ for some $\ell\in\mathbb{N}$, and that $x_{1}>n\varepsilon$. If the event $E_{n,\varepsilon}(\boldsymbol{p}^{(1)},\boldsymbol{p}^{(2)})$ holds, then any optimal alignment aligns all but all but at most an $\varepsilon$ proportion of the integer intervals $[v(k-1)+1,vk]$, $k=1,\ldots,d$, in $\mathbf{X}^{(1)}_{n}$, to integer intervals
$[r_{k-1}^{(j)}+1,r_{k}^{(j)}]$ in $\mathbf{X}^{(j)}_{n}$ of length no smaller than $vp_{j}^{(1)}$, for all $j=2,\ldots,m$. For each fixed $j\in\{2,\ldots,m\}$, the maximum number of integer intervals in $\mathbf{X}^{(1)}_{n}$ which could be matched with integer intervals in $\mathbf{X}^{(j)}_{n}$ of length less than $vp_{j}^{(1)}$ is at most $\varepsilon d$. Between $0$ and $x_{1}$, there are $\ell$ intervals from the
partition $[v(k-1)+1,vk]$, $k=1,\ldots,d$. Therefore, at least $\ell-\varepsilon d$ of these intervals are matched to integer intervals in $\mathbf{X}^{(j)}_{n}$ of length no less than $vp_{j}^{(1)}$. As a consequence, if the event $E_{n,\varepsilon}(\boldsymbol{p}^{(1)},\boldsymbol{p}^{(2)})$ holds, $x_{1}$ can get matched by an optimal alignment to a value of $x_{j}$ no less
than $(\ell-\varepsilon d)vp_{j}^{(1)}<p_{j}^{(1)}x_{1}-n\varepsilon p_{j}^{(1)}$. Finally, if $x_{1}$ is not an integer multiple of $v$, let $\tilde{x}_{1}$ be the largest integer multiple of $v$ which is smaller than $x_{1}$, then $x_{1}-\tilde{x}_{1}<v$. Note that, for each
$j\in\{2,\ldots,m\}$, $x_{1}$ get aligned with a point in the $j$-th sequence which cannot be before the point where $\tilde{x}_{1}$ gets aligned to in the same sequence. From the above discussion, since $\tilde{x}_{1}$ is an integer multiple of $v$, it gets aligned to a point in the $j$-the sequence which is no less than $p_{j}^{(1)}\tilde{x}_{1}-n\varepsilon p_{j}^{(1)}>p_{j}^{(1)}x_{1}-n\varepsilon p_{j}^{(1)}-p_{j}^{(1)}v$. Hence, if the event $E_{n,\varepsilon}(\boldsymbol{p}^{(1)},\boldsymbol{p}^{(2)})$ holds, then for each $j\in\{2,\ldots,m\}$, $x_{1}$ get aligned to a point
in the $j$-th sequence which is ``above or on" $p_{j}^{(1)}x_{1}-n\varepsilon p_{j}^{(1)}-p_{j}^{(1)}v$. This completes the proof of the first inclusion in \eqref{eq:InclEnDn}.

Therefore, by \eqref{eq:EstProbCompEn}, and for each $j=2,\ldots,m$,
\begin{align*}
\mathbb{P}\left(D_{1,j}^{c}\right)\leq\mathbb{P}\left(E_{n,\varepsilon}^{c}(\boldsymbol{p}^{(1)},\boldsymbol{p}^{(2)})\right)\leq\exp\left(-n\left(-\frac{1+\ln(1+v)}{v}(m-1)+\frac{\eta^{2}\varepsilon^{2}}{2m^{3}D^{2}}\right)\right).
\end{align*}
Similarly, the second inclusion in \eqref{eq:InclEnDn} implies that, for each $j=2,\ldots,m$,
\begin{align*}
\mathbb{P}\left(D_{2,j}^{c}\right)\leq\mathbb{P}\left(E_{n,\varepsilon}^{c}(\boldsymbol{p}^{(1)},\boldsymbol{p}^{(2)})\right)\leq\exp\left(-n\left(-\frac{1+\ln(1+v)}{v}(m-1)+\frac{\eta^{2}\varepsilon^{2}}{2m^{3}D^{2}}\right)\right).
\end{align*}
Finally, we note that
\begin{align*}
\mathbb{P}\left(D_{n,\varepsilon}^{c}(\boldsymbol{p}^{(1)},\boldsymbol{p}^{(2)})\right)\leq\sum_{j=2}^{m}\left(\mathbb{P}\left(D_{1,j}^{c}\right)+\mathbb{P}\left(D_{2,j}^{c}\right)\right),
\end{align*}
which completes the proof of the theorem.\hfill $\Box$
\begin{remark}\label{rem:CloseDiagMultiSeqs}
The proofs of Theorem \ref{thm:CloseDiagonal1} and Theorem \ref{thm:CloseDiagonal2} depend neither on the size of the alphabet nor on the choice of the score function. Hence, both theorems remain true for countably infinite alphabets, and for arbitrary permutation-invariant score functions as long as \eqref{eq:BdScore} and \eqref{eq:BddCondScore} are satisfied. For continuous (uncountable) alphabets, both theorems still remain true provided one has $\gamma^{*}>0$.
\end{remark}

\section{The Main Result}\label{sec:mainresult}

Recall that the Kolmogorov distance $d_{K}$, between two probability distributions $\nu_{1}$ and $\nu_{2}$ on $\mathbb{R}$, is defined as
\begin{align*}
d_{K}(\nu_{1},\nu_{2})=\sup_{h\in\mathcal{H}_{1}}\left|\int h\,d\nu_{1}-\int h\,d\nu_{2}\right|,
\end{align*}
where $\mathcal{H}_{1} :=\{\mathbf{1}_{(-\infty, x]}:\,x\in\mathbb{R}\}$.

Our main result is then:
\begin{theorem}\label{thm:CLT}
Let $\mathbf{X}^{(j)}_{n}=(X^{(j)}_{1},\ldots,X^{(j)}_{n})$, $j=1,\ldots,m$, be $m$ independent sequences of i.i.d. random variables drawn from a finite alphabet $\mathcal{A}\subset\mathbb{R}$. Let the score function $S$ satisfy \eqref{eq:BdScore} and \eqref{eq:BddCondScore} and let further
\begin{align}\label{eq:VarLowerBound}
\text{Var}(L_{n})\geq C_{*}n,\quad\text{for all }\,n\in\mathbb{N},
\end{align}
and some constant $C_{*}>0$, independent of $n$. Then, for all $n\in\mathbb{N}$,
\begin{align}\label{eq:BoundWasserstein}
d_{K}\left(\frac{L_{n}-\mathbb{E}(L_{n})}
{\sqrt{\text{Var}(L_{n})}},\,\mathcal{G}\right)\leq C\frac{(\ln n)^{3/4}}{n^{3/14}},
\end{align}
where $\mathcal{G}$ is a standard normal random variable, and where $C>0$ is a constant independent of $n$. In fact, if $\text{Var}(L_{n})\geq C_{*}(\ln n)^{\beta}n^{11/14}$ for some $\beta>3/4$, then a self-normalized CLT remains valid for $L_{n}$.
\end{theorem}
\begin{remark}\label{rem:VarLowerBound}
Below, we illustrate some instances where the variance lower bound \eqref{eq:VarLowerBound} holds true.
\begin{itemize}
\item [(i)] Consider the standard  LCS problem with $m=2$, i.e., let $(X_{i})_{i\geq 1}$ and $(Y_{i})_{i\geq 1}$ be two independent sequences of i.i.d. random variables with values in a finite alphabet $\mathcal{A}=\{\alpha_{1},\ldots,\alpha_{m}\}$, and let
    \begin{align*}
    \mathbb{P}\left(X_{1}=\alpha_{k}\right)=\mathbb{P}\left(Y_{1}=\alpha_{k}\right)=:p_{k},\quad k=1,\ldots,m.
    \end{align*}
    Next, let $p_{j_{0}}>1/2$, for some $j_{0}\in\{1,...,m\}$, $K=\min(2^{-4}10^{-2}e^{-67},1/800m)$, and let $\max_{j\neq j_{0}}p_{j}\leq\min\{2^{-2}e^{-5}K/m,K/2m^{2}\}$. Then, as shown in~\cite{HoudreMa:2014}, the variance of $L_{n}$ is of order $n$. The work  of~\cite{HoudreMa:2014} is based on the approach of~\cite{LemberMatzinger:2009} (and the references therein) where the order of the variance is also shown to be linear in the case of asymmetric Bernoulli random variables.
\item [(ii)] The linear order of the variance  is also proved in a score function setting with binary alphabets in~\cite{HoudreMatzinger:2007}. As noted in~\cite{HoudreMa:2014}, these results can be extended to multiple sequences on finite alphabets.
\item [(iii)] There are also instances where the variance is shown to be of order $n$ in a dependent random word setting, see e.g., \cite{AmsaluHoudreMatzinger:2012}. However, the proof below will not work for such cases, as the independence is crucial in the argument.
\end{itemize}
\end{remark}

The proof of Theorem \ref{thm:CLT} relies on a recent development in Stein's method due to Lachi\`eze-Rey and Peccati~\cite{LachiezeReyPeccati:2015} along with the results obtained in the previous sections (See~\cite{ChenGoldsteinShao:2011} and~\cite{Ross:2011} for thorough surveys on Stein's method).

To state the required results, some more notations are needed. Let $\boldsymbol{W}=(W_{1},\ldots,W_{n})$ and $\boldsymbol{W}'=(W_{1}',\ldots,W_{n}')$ be two i.i.d. $\mathbb{R}^{n}$-valued random vectors whose components are also independent. For $A\subset[n]$, define the random vector $\boldsymbol{W}^{A}$ by
\begin{align*}
W_{i}^{A}:=\left\{\begin{array}{ll} W_{i}' &\text{if }\,i\in A \\ W_{i} &\text{if }\,i\notin A \end{array}\right..
\end{align*}
For $A=\{j\}$, we denote $W^{j}$ instead of $W^{\{j\}}$ for further ease of notation.

For a given Borel measurable function $f:\mathbb{R}^{n}\rightarrow\mathbb{R}$ and $A \subset [n]$, let
\begin{align*}
T_{A}:=\sum_{j\notin A}\Delta_{j}f(\boldsymbol{W})\Delta_{j} f(\boldsymbol{W}^{A})\qquad\text{and}\qquad T_{A}':=\sum_{j \notin A}\Delta_{j} f(\boldsymbol{W})\left|\Delta_{j}f(\boldsymbol{W}^{A})\right|,
\end{align*}
where
\begin{align*}
\Delta_{j}f(\boldsymbol{W}):=f(\boldsymbol{W})-f(\boldsymbol{W}^{j}).
\end{align*}
Finally, let
\begin{align*}
T:=\frac{1}{2}\sum_{A\subsetneqq [n]}\frac{T_{A}}{\binom{n}{|A|}(n-|A|)}\qquad\text{and}\qquad T':=\frac{1}{2}\sum_{A\subsetneqq [n]}\frac{T_{A}'}{\binom{n}{|A|}(n-|A|)},
\end{align*}
where $|A|$ denotes the cardinality of $A$.
\begin{theorem}\label{thm:PR} (Lachi\`eze-Rey and Peccati~\cite[Theorem 4.2]{LachiezeReyPeccati:2015})$\,$
Let all the terms be defined as above, and let $0<\sigma^{2}:= \text{Var}(f(\boldsymbol{W}))<\infty$. Then,
\begin{align}
d_{K}\!\!\left(\!\frac{f(\boldsymbol{W})-\mathbb{E}(f(\boldsymbol{W}))}{\sqrt{\text{Var}(f(\boldsymbol{W})})},\,\mathcal{G}\!\right)&\leq\frac{1}{\sigma^{2}}\sqrt{\text{Var}\left(\mathbb{E}(T|\boldsymbol{W})\right)}+\frac{1}{\sigma^{2}}\sqrt{\text{Var}\left(\mathbb{E}(T'|\boldsymbol{W})\right)}\nonumber\\
\label{eq:PRbound} &\quad\,+\!\frac{1}{4\sigma^{3}}\!\sum_{j=1}^{n}\!\!\sqrt{\mathbb{E}\!\left(\left|\Delta_{j}f(\boldsymbol{W})\right|^{6}\right)}\!+\!\frac{\sqrt{2\pi}}{16\sigma^{3}}\!\sum_{j=1}^{n}\!\!\sqrt{\mathbb{E}\!\left(\left|\Delta_{j}f(\boldsymbol{W})\right|^{3}\right)},
\end{align}
where $\mathcal{G}$ is a standard normal random variable.
\end{theorem}

We conclude this section by noting that the result of Lachi\`eze-Rey and Peccati in Theorem \ref{thm:PR} is akin to the corresponding main result of  Chatterjee~\cite[Theorem 2.2]{Chatterjee:2008}, where a similar bound is obtained for the Wasserstein distance instead of the Kolmogorov one, and is therefore an improvement when directly estimating the Kolmogorov distance. The work in~\cite{HoudreIslak:2014} on the LCS problem used Chatterjee's result, and thus the results below also provide an improved normal convergence rate for $L_{n}$ in the LCS case. Note, however, that the major achievement of~\cite{HoudreIslak:2014} is in providing a first normal convergence result for $L_{n}$, rather than obtaining a rate of convergence result.

\noindent
\paragraph{Proof of Theorem \ref{thm:CLT}.} Let us make two comments before beginning the proof. First, in the proof, a constant $C$ may vary from an expression to another. Note, however, that $C$ will always be independent of $n$, the length of the sequences. Second, we do not worry about having quantities like $\ln n$, $n^{\alpha}$, etc., which should actually be $[n^{\alpha}]$, $[\ln n]$, etc.. This does not cause any problems as we are interested in asymptotic bounds.

We are now ready to begin the proof of Theorem \ref{thm:CLT}. To do so, let
\begin{align*}
\boldsymbol{W}=\left(\mathbf{X}_{n}^{(1)},\,\cdots,\mathbf{X}_{n}^{(m)}\right)\qquad\text{and}\qquad f(\boldsymbol{W})=L_{n}\left(\mathbf{X}_{n}^{(1)};\,\cdots;\mathbf{X}_{n}^{(m)}\right).
\end{align*}
We begin with estimating the last two terms on the right-hand side of \eqref{eq:PRbound}. Recalling \eqref{eq:VarLowerBound} and since $|\Delta_{j}f(\boldsymbol{W})|\leq s^{*}$ for any $j\in [n]$, we obtain that
\begin{align}\label{eq:EstLastTwoPR}
\frac{1}{4\sigma^{3}}\sum_{j=1}^{n}\sqrt{\mathbb{E}\left(\left|\Delta_{j}f(\boldsymbol{W})\right|^{6}\right)}+\frac{\sqrt{2\pi}}{16\sigma^{3}}\sum_{j=1}^{n}\sqrt{\mathbb{E}\left(\left|\Delta_{j}f(\boldsymbol{W})\right|^{3}\right)}\leq C\frac{1}{\sqrt n}.
\end{align}

Next, let us move to the estimation of the conditional variance terms in \eqref{eq:PRbound}, and make three more simple remarks:
\begin{itemize}
\item [(i)] Since $\text{Var}(\mathbb{E}(T|\boldsymbol{W}))\leq\text{Var}(T)$ and $\text{Var}(\mathbb{E}(T'|\boldsymbol{W}))\le\text{Var}(T')$, we may and do focus on estimating $\text{Var}(T)$ and $\text{Var}(T')$.
\item [(ii)] The estimate for $\text{Var}(T')$ is quite similar to the one for $\text{Var}(T)$. Henceforth, we will primarily work on $\text{Var}(T)$ and just briefly discuss the differences in the other case.
\item [(iii)] Parts of the estimation of $\text{Var}(T)$ follows along the lines of~\cite{HoudreIslak:2014}, and therefore the details will be skipped at certain steps.
\end{itemize}

Continuing the proof, set
\begin{align*}
\mathcal{S}_{1}:=\left\{(A,B,j,k):\,A\subsetneqq [2n],\,B\subsetneqq [2n],\,j\notin A,\,k \notin B\right\}.
\end{align*}
Introducing the notations
\begin{align*}
\kappa_{n,A}=\frac{1}{(n-|A|)\binom{n}{|A|}}\qquad\text{and}\qquad\kappa_{n,A,B}=\frac{1}{(n-|A|)\binom{n}{|A|}(n-|B|)\binom{n}{|B|}},
\end{align*}
$\text{Var}(T)$ can then be expressed as
\begin{align}
\text{Var}(T)&=\frac{1}{4}\text{Var}\left(\sum_{A\subsetneqq [mn]}\sum_{j\notin A}\kappa_{n,A}\Delta_{j}f( \boldsymbol{W})\Delta_{j}f(\boldsymbol{W}^{A})\right)\nonumber\\
\label{eq:VarExpanT} &=\frac{1}{4}\sum_{(A,B,j,k)\in\mathcal{S}_1}\kappa_{n,A,B}\,\text{Cov}\left(\Delta_{j}f(\boldsymbol{W})\Delta_{j}f(\boldsymbol{W}^{A}),\,\Delta_{k}f(\boldsymbol{W})\Delta_{k}f(\boldsymbol{W}^{B})\right).
\end{align}

Our strategy is now to further divide $\mathcal{S}_{1}$ into pieces and then to estimate the contributions of each piece separately. The following proposition, and a straightforward conditional version of it will be used repeatedly throughout. Its proof can be found in~\cite{HoudreIslak:2014}.
\begin{proposition}\label{prop:CombComp}
Let $\mathcal{R}$ be a subset of $[mn]^{2}$, and let
\begin{align*}
\mathcal{S}^{*}:=\left\{(A,B,j,k):\,A\subsetneqq [mn],\,B\subsetneqq [mn],\,j\notin A,\,k\notin B,\,(j,k)\in\mathcal{R}\right\}.
\end{align*}
Let $g:\mathcal{S}^{*}\rightarrow\mathbb{R}$ be such that $\|g\|_{\infty}\leq C$, then
\begin{align*}
\sum_{(A,B,j,k)\in\mathcal{S}^{*}}\left|\kappa_{mn,A,B}g(A,B,j,k)\right|\leq C\,|\mathcal{R}|.
\end{align*}
\end{proposition}

Taking $\mathcal{R}=[mn]^{2}$, Proposition \ref{prop:CombComp} yields the estimate
\begin{align*}
\sum_{(A,B,j,k)\in\mathcal{S}_{1}}\left(\kappa_{n,A,B}\,\text{Cov}\left(\Delta_{j}f(\boldsymbol{W})\Delta_{j}f(\boldsymbol{W}^{A}),\,\Delta_{k}f(\boldsymbol{W})\Delta_{k}f(\boldsymbol{W}^{B})\right)\right)\leq Cn^{2}.
\end{align*}
Hence, $\text{Var}(T)\leq Cn^{2}$ giving a suboptimal result for our purposes, and we therefore begin a detailed estimation study to improve the variance upper bound to $o(n^2)$.

Returning to the estimation of \eqref{eq:VarExpanT}, first, for notational convenience, we write $\sum_{1}$ in place of $\Sigma_{(A,B,j,k)\in\mathcal{S}_{1}}$ in the sequel. Also, for random variables $U$, $V$, and another random variable $Z$ taking its values in $R\subset\mathbb{R}$, and with another abuse of notation, we write
\begin{align*}
\text{Cov}_{Z=z}(U,V):=\mathbb{E}\left(\left(U-\mathbb{E}(U)\right)\left(V-\mathbb{E}(V)\right)|Z=z\right),\quad z\in R.
\end{align*}

Let, now, the random variable $Z$ be the indicator function of the event $E_{n,\varepsilon}(\boldsymbol{p}^{(1)},\boldsymbol{p}^{(2)})$, where $\varepsilon=c_{1}\sqrt{(1+\ln(1+n^{\alpha})/n^{\alpha}}$, with $c_{1}$ as in Remark \ref{rmk:epsilonchoice}. Then,
\begin{align}
&\sum_{1}\kappa_{mn,A,B}\,\text{Cov}\left(\Delta_{j}f(\boldsymbol{W})\Delta_{j}f(\boldsymbol{W}^{A}),\,\Delta_{k}f(\boldsymbol{W})\Delta_{k}f(\boldsymbol{W}^{B})\right)\nonumber\\
&\quad=\sum_{1}\kappa_{mn,A,B}\,\text{Cov}_{Z=0}\left(\Delta_{j}f(\boldsymbol{W})\Delta_{j}f(\boldsymbol{W}^{A}),\,\Delta_{k}f( \boldsymbol{W})\Delta_{k}f(\boldsymbol{W}^{B})\right)\mathbb{P}(Z=0)\nonumber\\
\label{eq:VarExpanTZ01} &\quad\quad\,+\sum_{1}\kappa_{mn,A,B}\,\text{Cov}_{Z=1}\left(\Delta_{j}f(\boldsymbol{W})\Delta_{j}f(\boldsymbol{W}^{A}),\,\Delta_{k}f(\boldsymbol{W})\Delta_{k}f(\boldsymbol{W}^{B})\right)\mathbb{P}(Z=1).
\end{align}
To estimate the first term on the right-hand side of \eqref{eq:VarExpanTZ01}, first note that
\begin{align*}
\text{Cov}_{Z=0}\left(\Delta_{j}f(\boldsymbol{W})\Delta_{j}f(\boldsymbol{W}^{A}),\,\Delta_{k}f(\boldsymbol{W})\Delta_{k}f(\boldsymbol{W}^{B})\right)\leq 4(s^{*})^{2},
\end{align*}
which, when combined with the estimate in \eqref{eq:EstProbCompEn} and Proposition \ref{prop:CombComp}, gives
\begin{align}
&\sum_{1}\kappa_{mn,A,B}\,\text{Cov}_{Z=0}\left(\Delta_{j}f(\boldsymbol{W})\Delta_{j}f(\boldsymbol{W}^{A}),\,\Delta_{k}f(\boldsymbol{W})\Delta_{k}f(\boldsymbol{W}^{B})\right)\mathbb{P}(Z=0)\nonumber\\
\label{eq:EstVarExpanTZ0} &\quad\leq 4(s^{*})^{2}n^{2}\,\exp\left(-n^{1-\alpha}(m-1)\left(1+\ln\left(1+n^{\alpha}\right)\right)\right).
\end{align}
For the second term on the right-hand side of \eqref{eq:VarExpanTZ01}, we begin with the trivial bound on $\mathbb{P}(Z=1)$ to obtain that
\begin{align}
&\sum_{1}\kappa_{mn,A,B}\,\text{Cov}_{Z=1}\left(\Delta_{j}f(\boldsymbol{W})\Delta_{j}f(\boldsymbol{W}^{A}),\,\Delta_{k}f(\boldsymbol{W})\Delta_{k}f(\boldsymbol{W}^{B})\right)\mathbb{P}(Z=1)\nonumber\\
\label{eq:EstVarExpanTZ11} &\quad\leq\sum_{1}\kappa_{mn,A,B}\,\text{Cov}_{Z=1}\left(\Delta_{j}f(\boldsymbol{W})\Delta_{j}f(\boldsymbol{W}^{A}),\,\Delta_{k}f(\boldsymbol{W})\Delta_{k}f(\boldsymbol{W}^{B})\right).
\end{align}

Finer decompositions are then needed to handle this last summation, and for this purpose, we specify an optimal alignment with certain properties. In the sequel, $\mathbf{r}$ denotes a uniquely defined optimal alignment which also specifies the $m$-tuples, in the sequences $\mathbf{X}^{(1)}_{n},\ldots,\mathbf{X}^{(m)}_{n}$, contributing to the optimal subsequences. Such an alignment always exists, and so we can define an injective map from $\mathbf{X}^{(1)}_{n},\ldots,\mathbf{X}^{(m)}_{n}$ to the set of alignments. This abstract construction is enough for our purposes, since the argument below is independent of the choice of the alignment. Another definition is then in order.
\begin{definition}\label{def:OptAlign}
For the optimal alignment $\mathbf{r}$, each of the sets
\begin{align*}
\left\{\mathbf{X}^{(1)}_{v(i-1)+1}\mathbf{X}^{(1)}_{v(i-1)+2}\cdots\mathbf{X}^{(1)}_{v(i)};\mathbf{X}^{(2)}_{r_{i-1}^{(2)}+1}\mathbf{X}^{(2)}_{r_{i-1}^{(2)}+2}\cdots\mathbf{X}^{(2)}_{r_{i}^{(2)}};\ldots;\mathbf{X}^{(m)}_{r_{i-1}^{(m)}+1}\mathbf{X}^{(m)}_{r_{i-1}^{(m)}+2}\cdots\mathbf{X}^{(m)}_{r_{i}^{(m)}}\right\},
\end{align*}
for $i=1,\ldots,d$, is called a \emph{cell} of $\mathbf{r}$.
\end{definition}

In particular, and clearly, any optimal alignment with $v=n^{\alpha}$ has $d=n^{1-\alpha}$ cells. Let us next introduce some more notation which will be used below. For any given $j\in [mn]$, let $P_{j}$ be the cell containing $W_{j}$ where, again, $\boldsymbol{W}=(\mathbf{X}_{n}^{(1)},\,\cdots,\mathbf{X}_{n}^{(m)})$.  Set $P_{j}=(P_{j}^{(1)};P_{j}^{(2)};\ldots;P_{j}^{(m)})$ where $P_{j}^{(i)}$ is the subword of $\mathbf{X}_{n}^{(i)}$ corresponding to $P_{j}$. Note that, for each $j\in [mn]$, $P_{j}^{(1)}$ contains $n^{\alpha}$ letters but that $P_{j}^{(i)}$, $i\geq 2$, might be empty. Here is an example.
\begin{example}
Take $n=12$, $m =3$ and $\mathcal{A}=[3]$. Let
\begin{align*}
\mathbf{X}^{(1)}_{n}=(1,1,2,1,2,1,1,2,1,1,3,1),\quad\text{and}\quad\mathbf{X}^{(2)}_{n}=\mathbf{X}^{(3)}_{n}=(2,1,1,3,2,3,1,2,1,1,1,1),
\end{align*}
and $ \boldsymbol{W}=(\mathbf{X}^{(1)}_{n},\mathbf{X}^{(2)}_{n},\mathbf{X}^{(3)}_{n})$. Then, $L_{12}=8$, and choosing $v=3$, the number of
cells in the optimal alignment is $d=4$. One possible choice for these cells is
$$(X_1^{(1)} X_2^{(1)} X_3^{(1)}; X_1^{(2)} X_2^{(2)} X_3^{(2)} X_4^{(2)} X_5^{(2)}; X_1^{(3)}
X_2^{(3)} X_3^{(3)} X_4^{(3)} X_5^{(3)}) =(1 1 2;2 1 1 3 2; 2 1 1 3 2),$$
$$(X_4^{(1)} X_5^{(1)} X_6^{(1)}; \emptyset; \emptyset) =(1 2 1; \emptyset; \emptyset),$$
$$(X_7^{(1)} X_8^{(1)} X_9^{(1)};  X_6^{(2)} X_7^{(2)} X_8^{(2)} X_9^{(2)}; X_6^{(3)} X_7^{(3)} X_8^{(3)} X_9^{(3)}) =(1 2 1;3 1 2 1; 3 1 2 1),$$ and
$$(X_{10}^{(1)} X_{11}^{(1)} X_{12}^{(1)}; X_{10}^{(2)} X_{11}^{(2)} X_{12}^{(2)}; X_{10}^{(3)} X_{11}^{(3)} X_{12}^{(3)}) =(1 3 1;1 1 1; 1 1 1).$$
For example, focusing on $W_8=X_8$, we have
\begin{align*}
P_{8}=\left(P_{8}^{(1)};P_{8}^{(2)};P_{8}^{(3)}\right)=(1 2 1; 3 1 2 1; 3 1 2 1).
\end{align*}
\end{example}

Returning to the proof of Theorem \ref{thm:CLT}, define the following subsets of $\mathcal{S}_{1}$ with respect to the alignment $\mathbf{r}$:
\begin{align*}
\mathcal{S}_{1,1}&=\left\{(A,B,j,k)\in\mathcal{S}_{1}:\,W_{j}\;\text{and}\;W_{k}\;\text{are in the same cell of}\;\mathbf{r}\right\},\\
\mathcal{S}_{1,2}&=\left\{(A,B,j,k)\in\mathcal{S}_{1}:\,W_{j}\;\text{and}\;W_{k}\;\text{are in different cells of}\;\mathbf{r}\right\}.
\end{align*}
Clearly, $\mathcal{S}_{1,1}\cap\mathcal{S}_{1,2}=\emptyset$ and $\mathcal{S}_{1}=\mathcal{S}_{1,1}\cup\mathcal{S}_{1,2}$. Now, for any $\mathcal{S}\subset\mathcal{S}_{1}$ and $(A,B,j,k)\in\mathcal{S}_{1}$, define
\begin{align*}
\text{Cov}_{Z=1,(A,B,j,k),\mathcal{S}}(X,Y):=\mathbb{E}\left(\left(X-\mathbb{E}(X)\right)\left(Y-\mathbb{E}(Y)\right){\bf 1}_{(A,B,j,k)\in\mathcal{S}}\,|\,Z=1\right).
\end{align*}
Also, write $\text{Cov}_{Z=1,\mathcal{S}}(X,Y)$ instead of $\text{Cov}_{Z=1,(A,B,j,k),\mathcal{S}}(X,Y)$ when the value of $(A,B,j,k)$ is clear from the context. The right-hand side of \eqref{eq:EstVarExpanTZ11} can then be further decomposed as
\begin{align}
&\sum_{1}\kappa_{mn,A,B}\,\text{Cov}_{Z=1}\left(\Delta_{j}f(\boldsymbol{W})\Delta_{j}f(\boldsymbol{W}^{A}),\,\Delta_{k}f(\boldsymbol{W})\Delta_{k}f(\boldsymbol{W}^{B})\right)\nonumber\\
&\quad=\sum_{1}\kappa_{mn,A,B}\,\text{Cov}_{Z=1,\mathcal{S}_{1,1}}\left(\Delta_{j}f(\boldsymbol{W})\Delta_{j}f(\boldsymbol{W}^{A}),\,\Delta_{k}f(\boldsymbol{W})\Delta_{k}f(\boldsymbol{W}^{B})\right)\nonumber\\
\label{eq:VarExpanTZ11S} &\quad\quad\,+\sum_{1}\kappa_{mn,A,B}\,\text{Cov}_{Z=1,\mathcal{S}_{1,2}}\left(\Delta_{j}f(\boldsymbol{W})\Delta_{j}f(\boldsymbol{W}^{A}),\,\Delta_{k}f(\boldsymbol{W})\Delta_{k}f(\boldsymbol{W}^{B})\right),
\end{align}
where, to further clarify the notation, note that, for example,
\begin{align*}
&\sum_{1}\kappa_{mn,A,B}\,\text{Cov}_{Z=1,\mathcal{S}_{1,1}}\left(\Delta_{j}f(\boldsymbol{W})\Delta_{j}f(\boldsymbol{W}^{A}),\,\Delta_{k}f(\boldsymbol{W})\Delta_{k}f(\boldsymbol{W}^{B})\right)\\
&\quad=\sum_{1}\mathbb{E}\left(\kappa_{mn,A,B}\,g(A,B,j,k){\bf 1}_{\{(A,B,j,k)\in\mathcal{S}_{1,1}\}}\,\Big|\,Z=1\right),
\end{align*}
where
\begin{align*}
g(A,B,j,k)&=\left(\Delta_{j}f(\boldsymbol{W})\Delta_{j}f(\boldsymbol{W}^{A})-\mathbb{E}\left(\Delta_{j}f(\boldsymbol{W})\Delta_{j}f(\boldsymbol{W}^{A})\right)\right)\\
&\quad\,\times\left(\Delta_{k}f(\boldsymbol{W})\Delta_{k}f(\boldsymbol{W}^{B})-\mathbb{E}\left(\Delta_{k}f(\boldsymbol{W})\Delta_{k}f(\boldsymbol{W}^{B})\right)\right).
\end{align*}

Let us now focus on the first term on the right-hand side of \eqref{eq:VarExpanTZ11S}. With $g$ as above,
\begin{align}
&\sum_{1}\kappa_{mn,A,B}\left|\text{Cov}_{Z=1,\mathcal{S}_{1,1}}\left(\Delta_{j}f(\boldsymbol{W})\Delta_{j}f(\boldsymbol{W}^{A}),\,\Delta_{k}f(\boldsymbol{W})\Delta_{k}f(\boldsymbol{W}^{B})\right)\right|\nonumber\\
&\quad\leq\mathbb{E}\left(\sum_{1}\kappa_{mn,A,B}\left|g(A,B,j,k)\right|{\bf 1}_{\{(A,B,j,k)\in\mathcal{S}_{1,1}\}}\,\Big|\,Z=1\right)\nonumber\\
\label{eq:EstVarExpanTZ11S111} &\quad\leq 4(s^{*})^{2}\mathbb{E}\left(\sum_{1}\kappa_{mn,A,B}\,{\bf 1}_{\{(A,B,j,k)\in\mathcal{S}_{1,1}\}}\,\Big|\,Z=1\right)=4(s^{*})^{2}\mathbb{E}\left(|\mathcal{R}||Z=1\right),
\end{align}
where
\begin{align*}
\mathcal{R}=\left\{(j,k)\in [mn]^{2}:\,W_{j}\;\text{and}\;W_{k}\;\text{are in the same cell of}\;\mathbf{r}\right\}.
\end{align*}
To estimate \eqref{eq:EstVarExpanTZ11S111}, for each $i=1,\ldots,d$, let $|\mathcal{R}_i|$ be the number of pairs of indices $(j,k)\in [mn]^{2}$ that are in the $i$-th cell, and let $\mathcal{G}_{i}$ be the event that $r_{i}^{(j)}-r_{i-1}^{(j)}\leq p_{j}^{(2)}n^{\alpha}$ for each $j=2,\ldots,m$. Then,
\begin{align}\label{eq:DecompRi}
\mathbb{E}\left(|\mathcal{R}|\,|\,Z=1\right)=\sum_{i=1}^{n^{1-\alpha}}\mathbb{E}\left(|\mathcal{R}_{i}|{\bf 1}_{\mathcal{G}_{i}}\,|\,Z=1\right)+\sum_{i=1}^{n^{1-\alpha}}\mathbb{E}\left(|\mathcal{R}_{i}|{\bf 1}_{\mathcal{G}_{i}^{c}}\,|\,Z=1\right).
\end{align}
For the first term on the right-hand side of \eqref{eq:DecompRi}, note that when $\mathcal{G}_{i}$ holds true, the $i$-th cell can contain at most $n^{\alpha}+(m-1)(\max_{j=2,\ldots,m}p_{j}^{(2)})n^{\alpha}=(1+(m-1)(\max_{j=2,\ldots,m}p_{j}^{(2)}))n^{\alpha}$ letters ($n^{\alpha}$ is for the letters in $\mathbf{X}^{(1)}_{n}$ and $(m-1)(\max_{j=2,\ldots,m}p_{j}^{(2)})n^{\alpha}$ is for the letters in the remaining sequences), and thus,
\begin{align*}
|\mathcal{R}_{i}|{\bf 1}_{\mathcal{G}_{i}}\leq\left(1+(m-1)\left(\max_{j=2,\ldots,m}p_{j}^{(2)}\right)\right)^{2}n^{2\alpha}.
\end{align*}
This gives
\begin{align}\label{eq:EstRi1}
\sum_{i=1}^{n^{1-\alpha}}\mathbb{E}\left(|\mathcal{R}_{i}|{\bf 1}_{\mathcal{G}_{i}}\,\big|\,Z=1\right)\leq\left(1+(m-1)\left(\max_{j=2,\ldots,m}p_{j}^{(2)}\right)\right)^{2}n^{1+\alpha}.
\end{align}

To estimate the second term on the right-hand side of \eqref{eq:DecompRi}, we begin with defining $K_{\varepsilon,\boldsymbol{p}^{(2)}}^n$ to be the event that the first cell is a subset of
\begin{align*}
\left(X^{(1)}_{1},\ldots,X^{(1)}_{n^{\alpha}};X^{(2)}_{1},\ldots,X^{(2)}_{(\max_{j}p^{(2)}_j)(n\varepsilon+2n^{\alpha})};\cdots;X^{(m)}_{1},\ldots,X^{(m)}_{(\max_{j}p^{(2)}_{j})(n\varepsilon+2n^{\alpha})}\right).
\end{align*}
Then,
\begin{align*}
E_{n,\varepsilon}\left(\boldsymbol{p}^{(1)},\boldsymbol{p}^{(2)}\right)\subset D_{n,\varepsilon}\left(\boldsymbol{p}^{(1)},\boldsymbol{p}^{(2)}\right)\subset K^{n}_{\varepsilon,\boldsymbol{p}^{(2)}}.
\end{align*}
Hence, we can now write
\begin{align*}
\sum_{i=1}^{n^{1-\alpha}}\mathbb{E}\left(|\mathcal{R}_{i}|{\bf 1}_{\mathcal{G}_{i}^{c}}|Z=1\right)=\sum_{i=1}^{n^{1-\alpha}}\mathbb{E}\left(|\mathcal{R}_{i}|{\bf 1}_{\mathcal{G}_{i}^{c}}{\bf 1}_{D_{n,\varepsilon}(\boldsymbol{p}^{(1)},\boldsymbol{p}^{(2)})}\Big|Z=1\right),
\end{align*}
and start by estimating the first element of the sum on the above right-hand side. Note that
\begin{align*}
|\mathcal{R}_{1}|{\bf 1}_{\mathcal{G}_{1}^{c}}{\bf 1}_{D_{n,\varepsilon}(\boldsymbol{p}^{(1)},\boldsymbol{p}^{(2)})}\leq |\mathcal{R}_{1}|{\bf 1}_{\mathcal{G}_{1}^{c}}{\bf 1}_{K^{n}_{\varepsilon}}\leq\left[n^{\alpha}+(m-1)\left(n\varepsilon\max_{j}p_{j}^{(2)}+2n^{\alpha}\max_{j}p_{j}^{(2)}\right)\right]^{2},
\end{align*}
which yields
\begin{align*}
\mathbb{E}\left(|\mathcal{R}_{1}|{\bf 1}_{\mathcal{G}_{1}^{c}}{\bf 1}_{D_{n,\varepsilon}(\boldsymbol{p}^{(1)},\boldsymbol{p}^{(2)})}\big|Z=1\right)\leq\left[\left(1+2(m-1)\max_{j}p_{j}^{(2)}\right)n^{\alpha}+(m-1)\max_{j}p_{j}^{(2)}n\varepsilon\right]^{2}.
\end{align*}
The remaining terms in the summation can be estimated in a similar way (see e.g., arguments in~\cite{HoudreIslak:2014}). Therefore, choosing $\varepsilon=(c_{1}^{2}(1+\ln(n^{\alpha}+1))/n^{\alpha})^{1/2}$ leads to
\begin{align*}
\sum_{i=1}^{n^{1-\alpha}}\mathbb{E}\left(|\mathcal{R}_{i}|{\bf 1}_{\mathcal{G}_{i}^{c}}\,|\,Z=1\right)\leq C_{1}n^{1+\alpha/2}\left(\ln n^{\alpha}\right)^{1/2}+C_{2}n^{2-\alpha}\ln n^{\alpha}+C_{3}n^{3-5\alpha/2}\left(\ln n^{\alpha}\right)^{3/2},
\end{align*}
where $C_{1}$, $C_{2}$ and $C_{3}$ are constants independent of $n$. Combining this last inequality with \eqref{eq:EstRi1} leads to \eqref{eq:DecompRi},
\begin{align}\label{eq:EstR}
\mathbb{E}\left(|\mathcal{R}|\,|\,Z=1\right)\leq C\left(n^{1+\alpha}+n^{1+\alpha/2}\left(\ln n^{\alpha}\right)^{1/2}+n^{2-\alpha}\ln n^{\alpha}+n^{3-5\alpha/2}\left(\ln n^{\alpha}\right)^{3/2}\right).
\end{align}
In turn, together with \eqref{eq:EstVarExpanTZ11S111}, \eqref{eq:EstR} yields the following estimate on the first sum on the right-hand side of \eqref{eq:VarExpanTZ11S}:
\begin{align}
&\sum_{1}\kappa_{mn,A,B}\left|\text{Cov}_{Z=1,\mathcal{S}_{1,1}}\left(\Delta_{j}f(\boldsymbol{W})\Delta_{j}f(\boldsymbol{W}^{A}),\,\Delta_{k}f(\boldsymbol{W})\Delta_{k}f(\boldsymbol{W}^{B})\right)\right|\nonumber\\
\label{eq:EstVarExpanTZ11S112} &\quad\leq C\left(n^{1+\alpha}+n^{1+\alpha/2}\left(\ln n^{\alpha}\right)^{1/2}+n^{2-\alpha}\ln n^{\alpha}+n^{3-5\alpha/2}\left(\ln n^{\alpha}\right)^{3/2}\right).
\end{align}

In order to estimate the second sum on the right-hand side of \eqref{eq:VarExpanTZ11S}, we need to decompose the covariance terms in such a way that the (conditional) independence of certain random variables occurs, simplifying the estimates themselves. For this purpose, for each $i\in [mn]$, let $f(P_{i})=L(P_{i})$ be the optimal score of $P_{i}^{(1)},\ldots,P_{i}^{(m)}$, and let $\widetilde{\Delta}_{i}f(\boldsymbol{W}):=f(P_{i})-f(P_{i}')$, where $P_{i}'$ is the same as $P_{i}$ except that $W_{i}$ is now replaced by the independent copy $W_{i}'$. Now for $(A,B,j,k)\in\mathcal{S}_{1}$,
\begin{align}
&\text{Cov}_{Z=1,\mathcal{S}_{1,2}}\left(\Delta_{j}f(\boldsymbol{W})\Delta_{j}f(\boldsymbol{W}^{A}),\,\Delta_{k}f(\boldsymbol{W})\Delta_{k}f(\boldsymbol{W}^{B})\right)\nonumber\\
&\quad=\text{Cov}_{Z=1,\mathcal{S}_{1,2}}\left(\left(\Delta_{j}f(\boldsymbol{W})-\widetilde{\Delta}_{j}f(\boldsymbol{W})\right)\Delta_{j}f(\boldsymbol{W}^{A}),\,\Delta_{k}f(\boldsymbol{W})\Delta_{k}f(\boldsymbol{W}^{B})\right)\nonumber\\
&\quad\quad\,+\text{Cov}_{Z=1,\mathcal{S}_{1,2}}\left(\widetilde{\Delta}_{j}f(\boldsymbol{W})\left(\Delta_{j}f(\boldsymbol{W}^{A})-\widetilde{\widetilde{\Delta}}_{j}f(\boldsymbol{W}^{A})\right),\,\Delta_{k}f(\boldsymbol{W})\Delta_{k}f(\boldsymbol{W}^{B})\right)\nonumber\\
&\quad\quad\,+\text{Cov}_{Z=1,\mathcal{S}_{1,2}}\left(\widetilde{\Delta}_{j}f(\boldsymbol{W})\widetilde{\widetilde{\Delta}}_{j}f(\boldsymbol{W}^{A}),\,\left(\Delta_{k}f(\boldsymbol{W})-\widetilde{\Delta}_{k}f(\boldsymbol{W})\right)\Delta_{k}f(\boldsymbol{W}^{B})\right)\nonumber\\
&\quad\quad\,+\text{Cov}_{Z=1,\mathcal{S}_{1,2}}\left(\widetilde{\Delta}_{j}f(\boldsymbol{W})\widetilde{\widetilde{\Delta}}_{j}f(\boldsymbol{W}^{A}),\,\widetilde{\Delta}_{k}f(\boldsymbol{W})\left(\Delta_{k}f(\boldsymbol{W}^{B})-\widetilde{\widetilde{\Delta}}_{k}f(\boldsymbol{W}^{B})\right)\right)\nonumber\\
\label{eq:DecompCovZ1S12} &\quad\quad\,+\text{Cov}_{Z=1,\mathcal{S}_{1,2}}\left(\widetilde{\Delta}_{j}f(\boldsymbol{W})\widetilde{\widetilde{\Delta}}_{j}f(\boldsymbol{W}^{A}),\,\widetilde{\Delta}_{k}f(\boldsymbol{W})\widetilde{\widetilde{\Delta}}_{k}f(\boldsymbol{W}^{B})\right),
\end{align}
where, for any $i\notin A$,
\begin{align*}
\widetilde{\widetilde{\Delta}}_{i}f(\boldsymbol{W}^{A})=f\left(\left.\boldsymbol{W}^{A}\right|_{P_{i}}\right)-f\left(\left.\boldsymbol{W}^{A\cup\{i\}}\right|_{P_{i}}\right),
\end{align*}
with $W^{A}|_{P_{i}}$ being the restriction of $W^{A}$ to the cell $P_{i}$. Above, we used the bilinearity of $\text{Cov}_{Z=1,\mathcal{S}_{1,2}}$ to express the left-hand side as a telescoping sum.

Let us begin by focusing on the last term on the right-hand side of \eqref{eq:DecompCovZ1S12}. Proceeding along the same lines as in~\cite{HoudreIslak:2014}, one shows that
\begin{align*}
&\left|\text{Cov}_{Z=1,\mathcal{S}_{1,2}}\left(\widetilde{\Delta}_{j}f(\boldsymbol{W})\widetilde{\widetilde{\Delta}}_{j}f(\boldsymbol{W}^{A}),\widetilde{\Delta}_{k}f(\boldsymbol{W})\widetilde{\widetilde{\Delta}}_{k}f(\boldsymbol{W}^{B})\right)\right|\\
&\quad\leq C\left(\mathbb{P}((A,B,j,k)\in\mathcal{S}_{1,1},Z=1)^{2}+\mathbb{P}(Z=0)\right).
\end{align*}
Hence,
\begin{align}
&\sum_{1}\kappa_{mn,A,B}\left|\text{Cov}_{Z=1,\mathcal{S}_{1,2}}\left(\widetilde{\Delta}_{j}f(\boldsymbol{W})\widetilde{\widetilde{\Delta}}_{j}f(\boldsymbol{W}^{A}),\,\widetilde{\Delta}_{k}f(\boldsymbol{W})\widetilde{\widetilde{\Delta}}_{k}
f(\boldsymbol{W}^{B})\right)\right|\nonumber\\
&\quad\leq C\sum_{1}\left(\kappa_{mn,A,B}\,\mathbb{P}\left((A,B,j,k)\in\mathcal{S}_{1,1},Z=1\right)^{2}+\kappa_{mn,A,B}\,\mathbb{P}(Z=0)\right)\nonumber\\
&\quad\leq C\sum_{1}\left(\kappa_{mn,A,B}\,\mathbb{P}\left((A,B,j,k)\in\mathcal{S}_{1,1},Z=1\right)+\kappa_{mn,A,B}\,\mathbb{P}(Z=0)\right)\nonumber\\
&\quad\leq C\sum_{1}\left(\kappa_{mn,A,B}\,\mathbb{P}\left((A,B,j,k)\in\mathcal{S}_{1,1}|Z=1\right)\mathbb{P}(Z=1)+\kappa_{mn,A,B}\,\mathbb{P}(Z=0)\right)\nonumber\\
\label{eq:EstVarExpanTZ11S12Last}
&\quad\leq C\!\left(n^{1+\alpha}\!+\!n^{1+\alpha/2}\sqrt{\ln n^{\alpha}}\!+\!n^{2-\alpha}\!\ln n^{\alpha}\!+\!n^{3-5\alpha/2}\!\left(\ln n^{\alpha}\right)^{3/2}\!\!+\!e^{-n^{1-\alpha}(m-1)\left(1+\ln\left(1+n^{\alpha}\right)\right)}\right),
\end{align}
where for the last step we made use of the estimates \eqref{eq:Zequalszeroestimate} and \eqref{eq:EstVarExpanTZ11S112}.

Next, we upper-bound the first summand in \eqref{eq:DecompCovZ1S12}. The other three terms can be dealt with in a similar way, and the details are thus omitted. Letting
\begin{align*}
U:=\left(\Delta_{j}f(\boldsymbol{W})-\widetilde{\Delta}_{j}f(\boldsymbol{W})\right)\Delta_{j}f(\boldsymbol{W}^{A}),\quad V:=\Delta_{k}f(\boldsymbol{W})\Delta_{k}f(\boldsymbol{W}^{B}),
\end{align*}
we proceed to estimating $\text{Cov}_{Z=1,\mathcal{S}_{1,2}}(U,V)$ by writing
\begin{align*}
\left|\text{Cov}_{Z=1,\mathcal{S}_{1,2}}(U,V)\right|&=\left|\mathbb{E}\left(\left(U-\mathbb{E}(U)\right)\left(V-\mathbb{E}(V)\right){\bf 1}_{(A,B,j,k)\in\mathcal{S}_{1,2}}|Z=1\right)\right|\\
&\leq\mathbb{E}\left(|UV|{\bf 1}_{\{(A,B,j,k)\in\mathcal{S}_{1,2}\}}|Z=1\right)+\mathbb{E}(|V|)\,\mathbb{E}\left(|U|{\bf 1}_{\{(A,B,j,k)\in\mathcal{S}_{1,2}\}}|Z=1\right)\\
&\quad\,+\mathbb{E}(|U|)\,\mathbb{E}\left(|V|{\bf 1}_{\{(A,B,j,k)\in\mathcal{S}_{1,2}\}}\,|\,Z=1\right)\\
&\quad\,+\mathbb{E}(|U|)\mathbb{E}(|V|)\,\mathbb{E}\left({\bf 1}_{\{(A,B,j,k)\in\mathcal{S}_{1,2}\}}\,|\,Z=1\right)\\
&=:T_{1}+T_{2}+T_{3}+T_{4}.
\end{align*}
Starting with the estimation of  $T_{1}$, since $|\Delta_{j}f( \boldsymbol{W}^{A})(\Delta_{k}f(\boldsymbol{W})\Delta_{k}f( \boldsymbol{W}^{B}))|\leq C$, we have
\begin{align}\label{eq:EstT1}
T_{1}\leq C\,\mathbb{E}_{Z=1}\left(\left|\Delta_{j}f(\boldsymbol{W})-\widetilde{\Delta}_{j}f(\boldsymbol{W})\right|{\bf 1}_{\{(A,B,j,k)\in\mathcal{S}_{1,2}\}}\right).
\end{align}
A similar estimate also reveals that
\begin{align}\label{eq:EstT2}
T_{2}\leq C\,\mathbb{E}_{Z=1}\left(\left|\Delta_{j}f(\boldsymbol{W})-\widetilde{\Delta}_{j}f(\boldsymbol{W})\right|{\bf 1}_{\{(A,B,j,k)\in\mathcal{S}_{1,2}\}}\right).
\end{align}
Next, for $T_{3}$ and $T_{4}$, by \eqref{eq:Zequalszeroestimate}, we have
\begin{align}
T_{3}+T_{4}\leq C\,\mathbb{E}(|U|)&\leq C\,\mathbb{E}\left(\left|\Delta_{j}f(\boldsymbol{W})-\widetilde{\Delta}_{j}f(\boldsymbol{W})\right|\right)\nonumber\\
&\leq C\,\mathbb{E}_{Z=1}\left(\left|\Delta_{j}f(\boldsymbol{W})-\widetilde{\Delta}_{j}f(\boldsymbol{W})\right|{\bf 1}_{\{(A,B,j,k)\in\mathcal{S}_{1,2}\}}\right)\nonumber\\
&\quad\,+C\,\mathbb{E}_{Z=1}\left(\left|\Delta_{j}f(\boldsymbol{W})-\widetilde{\Delta}_{j}f(\boldsymbol{W})\right|{\bf 1}_{\{(A,B,j,k)\in\mathcal{S}_{1,1}\}}\right)\nonumber\\
\label{eq:EstT3T4} &\quad\,+C\,e^{-n^{1-\alpha}(m-1)(1+\ln(1+n^{\alpha}))}.
\end{align}
Let $h(A,B,j,k)$ be the sum of the first four terms on the right-hand side of \eqref{eq:DecompCovZ1S12}. Then, performing similar estimations as in
getting \eqref{eq:EstT1}, \eqref{eq:EstT2} and \eqref{eq:EstT3T4}, for the second to fourth term of this sum, and observing that $|\Delta_{j}f(\boldsymbol{W})-\widetilde{\Delta}_{j}f(\boldsymbol{W})|$ and $|\Delta_{j}f(\boldsymbol{W}^{A})-\widetilde{\widetilde{\Delta}}_{j}f(\boldsymbol{W}^{A})|$,
$|\Delta_{k}f(\boldsymbol{W})-\widetilde{\Delta}_{k}f(\boldsymbol{W})|$ and $|\Delta_{k}f(\boldsymbol{W}^{B})-\widetilde{\widetilde{\Delta}}_{k}f(\boldsymbol{W}^{B})|$ are identically distributed, respectively, we obtain that
\begin{align*}
\sum_{1}\left|\kappa_{mn,A,B}\,h(A,B,j,k)\right|&\leq C\sum_{1}\kappa_{mn,A,B}\,\mathbb{E}_{Z=1}\left(\left|\Delta_{j}f(\boldsymbol{W})-\widetilde{\Delta}_{j}f(\boldsymbol{W})\right|{\bf 1}_{\{(A,B,j,k)\in\mathcal{S}_{1,2}\}}\right)\\
&\quad\,+C\sum_{1}\kappa_{mn,A,B}\,\mathbb{E}_{Z=1}\left(\left|\Delta_{j}f(\boldsymbol{W})-\widetilde{\Delta}_{j}f(\boldsymbol{W})\right|{\bf 1}_{\{(A,B,j,k)\in\mathcal{S}_{1,1}\}}\right)\\
&\quad\,+C\sum_{1}\kappa_{mn,A,B}\,\mathbb{E}_{Z=1}\left(\left|\Delta_{k}f(\boldsymbol{W})-\widetilde{\Delta}_{k}f(\boldsymbol{W})\right|{\bf 1}_{\{(A,B,j,k)\in\mathcal{S}_{1,2}\}}\right)\\
&\quad\,+C\sum_{1}\kappa_{mn,A,B}\,\mathbb{E}_{Z=1}\left(\left|\Delta_{k}f(\boldsymbol{W})-\widetilde{\Delta}_{k}f(\boldsymbol{W})\right|{\bf 1}_{\{(A,B,j,k)\in\mathcal{S}_{1,1}\}}\right)\\
&\quad\,+C\sum_{1}\kappa_{mn,A,B}\,e^{-n^{1-\alpha}(m-1)(1+\ln(1+n^{\alpha}))}.
\end{align*}
By making use of a symmetry argument, this gives
\begin{align*}
\sum_{1}\left|\kappa_{mn,A,B}\,h(A,B,j,k)\right|&\leq C\sum_{1}\kappa_{mn,A,B}\,\mathbb{E}_{Z=1}\left(\left|\Delta_{j}f(\boldsymbol{W})-\widetilde{\Delta}_{j}f(\boldsymbol{W})\right|{\bf 1}_{\{(A,B,j,k)\in\mathcal{S}_{1,2}\}}\right)\\
&\quad\,+C\sum_{1}\kappa_{mn,A,B}\,\mathbb{E}_{Z=1}\left(\left|\Delta_{j}f(\boldsymbol{W})-\widetilde{\Delta}_{j}f(\boldsymbol{W})\right|{\bf 1}_{\{(A,B,j,k)\in\mathcal{S}_{1,1}\}}\right)\\
&\quad\,+C\sum_{1}\kappa_{mn,A,B}\,e^{-n^{1-\alpha}(m-1)(1+\ln(1+n^{\alpha}))}.
\end{align*}
By Proposition \ref{prop:CombComp}, the third sum on the above right-hand side is upper-bounded by
\begin{align}\label{eq:EstCovZ1S123}
Cn^{2}\,e^{-n^{1-\alpha}(m-1)(1+\ln(1+n^{\alpha}))}.
\end{align}
Moreover, by \eqref{eq:EstR} and since $|\Delta_{j}f(\boldsymbol{W})-\widetilde{\Delta}_{j}f(\boldsymbol{W})|\leq C$, the middle sum, above, is
itself upper-bounded by
\begin{align}\label{eq:EstCovZ1S122}
C\left(n^{1+\alpha}+n^{1+\alpha/2}\left(\ln n^{\alpha}\right)^{1/2}+n^{2-\alpha}\ln n^{\alpha}+n^{3-5\alpha/2}\left(\ln n^{\alpha}\right)^{3/2}\right).
\end{align}
Finally, still following the proof in~\cite{HoudreIslak:2014}, the first sum, above, can be upper-bounded by
\begin{equation}\label{eq:EstCovZ1S121}
C\,n^{2}e^{-n^{1-\alpha}(m-1)(1+\ln(1+n^{\alpha}))},
\end{equation}
noting that~\cite[Proposition 2.2]{HoudreIslak:2014} continues to hold for the score function $S$ satisfying \eqref{eq:BddCondScore} and \eqref{eq:bddconditionK}.

Combining \eqref{eq:EstVarExpanTZ0}, \eqref{eq:EstVarExpanTZ11S112}, \eqref{eq:EstVarExpanTZ11S12Last}, \eqref{eq:EstCovZ1S123}, \eqref{eq:EstCovZ1S122} and \eqref{eq:EstCovZ1S121}, gives
\begin{align}
\text{Var}(T)&\leq C\left(n^{1+\alpha}+n^{1+\alpha/2}\left(\ln n^{\alpha}\right)^{1/2}+n^{2-\alpha}\ln n^{\alpha}+n^{3-5\alpha/2}\left(\ln n^{\alpha}\right)^{3/2}\right.\nonumber\\
\label{eq:EstVarT} &\quad\qquad +n^{2}\,e^{-n^{1-\alpha}(m-1)(1+\ln(1+n^{\alpha}))}\Big).
\end{align}

The treatment of the estimation of $\text{Var}(T')$ is very similar, and thus we only include a sketch. Again we first write
\begin{align}
\text{Var}(T')&=\frac{1}{4}\,\text{Var}\left(\sum_{A\subsetneqq [mn]}\sum_{j\notin A}\kappa_{n,A}\,\Delta_{j}f(\boldsymbol{W})\left|\Delta_{j}f(\boldsymbol{W}^{A})\right|\right)\nonumber\\
&=\frac{1}{4}\sum_{(A,B,j,k)\in\mathcal{S}_{1}}\kappa_{n,A,B}\,\text{Cov}\left(\Delta_{j}f(\boldsymbol{W})\left|\Delta_{j}f(\boldsymbol{W}^{A})\right|,\,\Delta_{k}f(\boldsymbol{W})\left|\Delta_{k}f(\boldsymbol{W}^{B})\right|\right)\nonumber\\
&=\frac{1}{4}\sum_{(A,B,j,k)\in\mathcal{S}_{1,1}}\kappa_{n,A,B}\,\text{Cov}\left(\Delta_{j}f(\boldsymbol{W})\left|\Delta_{j}f(\boldsymbol{W}^{A})\right|,\,\Delta_{k}f(\boldsymbol{W})\left|\Delta_{k}f(\boldsymbol{W}^{B})\right|\right)\nonumber\\
\label{eq:DecompVarT'} &\quad\,+\frac{1}{4}\sum_{(A,B,j,k)\in\mathcal{S}_{1,2}}\!\!\kappa_{n,A,B}\,\text{Cov}\left(\Delta_{j}f(\boldsymbol{W})\left|\Delta_{j}f(\boldsymbol{W}^{A})\right|,\,\Delta_{k}f(\boldsymbol{W})\left|\Delta_{k}f(\boldsymbol{W}^{B})\right|\right).
\end{align}
The estimation of the first summation of \eqref{eq:DecompVarT'} is exactly the same as its counterpart in  estimating $\text{Var}(T)$. For the estimation of the second term of \eqref{eq:DecompVarT'}, the following decomposition, which is analogous to \eqref{eq:DecompCovZ1S12}, is needed:
\begin{align}
&\text{Cov}_{Z=1,\mathcal{S}_{1,2}}\left(\Delta_{j}f(\boldsymbol{W})\left|\Delta_{j}f(\boldsymbol{W}^{A})\right|,\,\Delta_{k}f(\boldsymbol{W})\left|\Delta_{k}f(\boldsymbol{W}^{B})\right|\right)\nonumber\\
&\quad=\text{Cov}_{Z=1,\mathcal{S}_{1,2}}\left(\left(\Delta_{j}f(\boldsymbol{W})-\widetilde{\Delta}_{j}f(\boldsymbol{W})\right)\left|\Delta_{j}f(\boldsymbol{W}^{A})\right|,\,\Delta_{k}f(\boldsymbol{W})\left|\Delta_{k}f(\boldsymbol{W}^{B})\right|\right)\nonumber\\
&\quad\quad\,+\text{Cov}_{Z=1,\mathcal{S}_{1,2}}\left(\widetilde{\Delta}_{j}f(\boldsymbol{W})\left(\left|\Delta_{j}f(\boldsymbol{W}^{A})\right|-\left|\widetilde{\widetilde{\Delta}}_{j}f(\boldsymbol{W}^{A})\right|\right),\,\Delta_{k}f(\boldsymbol{W})\left|\Delta_{k}f(\boldsymbol{W}^{B})\right|\right)\nonumber\\
&\quad\quad\,+\text{Cov}_{Z=1,\mathcal{S}_{1,2}}\left(\widetilde{\Delta}_{j}f(\boldsymbol{W})\left|\widetilde{\widetilde{\Delta}}_{j}f(\boldsymbol{W}^{A})\right|,\,\left(\Delta_{k}f(\boldsymbol{W})-\widetilde{\Delta}_{k}f(\boldsymbol{W})\right)\left|\Delta_{k}f(\boldsymbol{W}^{B})\right|\right)\nonumber\\
&\quad\quad\,+\text{Cov}_{Z=1,\mathcal{S}_{1,2}}\left(\widetilde{\Delta}_{j}f(\boldsymbol{W})\left|\widetilde{\widetilde{\Delta}}_{j}f(\boldsymbol{W}^{A})\right|,\,\widetilde{\Delta}_{k}f(\boldsymbol{W})\left(\left|\Delta_{k}f(\boldsymbol{W}^{B})\right|-\left|\widetilde{\widetilde{\Delta}}_{k}f(\boldsymbol{W}^{B})\right|\right)\right)\nonumber\\
\label{eq:DecompT'CovZ1S12} &\quad\quad\,+\text{Cov}_{Z=1,\mathcal{S}_{1,2}}\left(\widetilde{\Delta}_{j}f(\boldsymbol{W})\left|\widetilde{\widetilde{\Delta}}_{j}f(\boldsymbol{W}^{A})\right|,\,\widetilde{\Delta}_{k}f(\boldsymbol{W})\left|\widetilde{\widetilde{\Delta}}_{k}f(\boldsymbol{W}^{B})\right|\right),
\end{align}
All the summands in \eqref{eq:DecompT'CovZ1S12} can be handled with arguments similar to the ones used for estimating the counterparts in $\text{Var}(T)$. The only difference is the use of certain elementary inequalities such as $\big||\Delta_{j}f(\boldsymbol{W}^{A})|-|\widetilde{\widetilde{\Delta}}_{j}f(\boldsymbol{W}^{A})|\big|\leq|\Delta_{j}f(\boldsymbol{W}^{A})-\widetilde{\widetilde{\Delta}}_{j}f(\boldsymbol{W}^{A})|$.
Combining all the estimates again leads to
\begin{align}
\text{Var}(T')&\leq C\left(n^{1+\alpha}+n^{1+\alpha/2}\left(\ln n^{\alpha}\right)^{1/2}+n^{2-\alpha}\ln n^{\alpha}+n^{3-5\alpha/2}\left(\ln n^{\alpha}\right)^{3/2}\right.\nonumber\\
\label{eq:EstVarT'} &\quad\qquad +n^{2}\,e^{-n^{1-\alpha}(m-1)(1+\ln(1+n^{\alpha}))}\Big).
\end{align}

Therefore, \eqref{eq:BoundWasserstein} follows immediately from Theorem \ref{thm:PR}, \eqref{eq:EstLastTwoPR}, \eqref{eq:EstVarT} and \eqref{eq:EstVarT'} with the choice $\alpha=4/7$.\hfill $\square$

\begin{remark}\label{rem:MainThm}
\begin{itemize}
\item [(i)] The constant $C$ in Theorem~\ref{thm:CLT} is independent of $n$, but depends on $m$, $s^{*}$ and $(\boldsymbol{p}^{(1)},\boldsymbol{p}^{(2)})$.
\item [(ii)] As already noted, the conclusion of Theorem~\ref{thm:CLT}, in the special case of the LCS problem, improves, to $(\ln n)^{3/4}/n^{3/14}$, the corresponding convergence rate obtained in~\cite{HoudreIslak:2014}. (The estimate on $d_{K}((L_{n}-\mathbb{E}(L_n))/\sqrt{\text{Var}(L_{n})},\mathcal{G})$ obtained there is of order $(\ln n)^{3/8}/n^{3/28}$.)
\item [(iii)] Of course, there is no reason for our rate to be sharp. Already, instead of the choice $v=n^{\alpha}$, a choice such as $v=h(n)$, for some optimal function $h$ would improve the rate.
\end{itemize}
\end{remark}

\section{Concluding remarks}\label{sec:Conclusion}

\begin{enumerate}
\item \emph{Generalized Bernoulli matching.} Note that for a given score function $S$, as in the case of the longest common subsequence problem, one can find the value of $L_{n}$ using the following recursion:
    \begin{align*}
    L_{i,j}=\max\left\{L_{i-1,j},L_{i,j-1},L_{i-1,j-1}+S(X_{i},Y_{j})\right\},\quad L(0,0)=L(i,0)=L(0,j)=0,
    \end{align*} 
    and recording the value $L_{n,n}$ as $L_{n}$. Note also that the score terms $S(X_{i},Y_{j})$ in the above recursion are dependent. In~\cite{MajumdarNechaev:2005}, Majumdar and Nechaev, by focusing on the finite discrete alphabet case and using the standard score function $S(a,b)={\bf 1}_{\{a=b\}}$, show that the Tracy-Widom law is the limiting distribution when the score terms of the recursion are assumed to be independent. The model considered by Majumdar and Nechaev is known as the \emph{Bernoulli matching model}, and it is natural to wonder whether or not their result remains true for a general scoring function.
\item \emph{Score functions over permutations.} It is shown in~\cite{HoudreIslak:2014} that the length of the longest common subsequences of two independent and uniform random permutations of $[n]$ converges in law to the Tracy-Widom distribution. It is thus again natural to question whether or not this result extends to a broader class of score functions. To focus on a specific example, let $\pi^{(n)}$ and $\rho^{(n)}$ be two independent uniform random permutations of $[n]$, and let the score function be $S(a,b)=1-\frac{|a-b|}{n-1}{\bf 1}_{\{|a-b|\leq c\}}$ for some real $c\geq 0$ which may or may not depend on $n$. Then, 
    \begin{align*}
    L_{n}=\max_{\boldsymbol{\mu},\boldsymbol{\nu}}\sum_{i=1}^{k}\left(1-\frac{|\pi^{(n)}_{\mu_{i}}-\rho^{(n)}_{\nu_{i}}|}{n-1}{\bf 1}_{\left\{\left|\pi^{(n)}_{\mu_{i}}-\rho^{(n)}_{\nu_{i}}\right|\leq c\right\}}\right).
    \end{align*}
    Now, for $c=0$ celebrated results of Vershik and Kerov~\cite{VershikKerov:1977} and Baik, Deift and Johansson~\cite{BaikDeiftJohansson:1999} respectively imply that the expectation of $L_{n}$ is of order $\sqrt{n}$ while the limiting law is the Tracy-Widom distribution. (We refer to Romik's book~\cite{Romik:2014} for an up to date account of this problem with a complete bibliography.) It is natural to conjecture that the same result will hold true when $c$ is any positive constant. However, as the value of $c$ grows with $n$ one would expect a transition in both the order of expectation and the limiting law. For instance, when $c=n-1$ it can easily be shown that the expectation is of order $n$ and it would not be surprising to have the Gaussian distribution as the limiting one.
\item \emph{Dependent setting.} Another problem of interest is the extension of the above results to the cases where dependencies do exist among the sequences $\mathbf{X}^{(i)}_{n}=(X^{(i)}_{1},\ldots,X^{(i)}_{n})$, $i=1,\ldots,m$. In particular, it is worth exploring whether or not the central limit theorem continues to hold when $\mathbf{X}^{(i)}_{n}$, $i=1,\ldots,m$ are Markov chains that are independent of each other.
\end{enumerate}

\end{document}